%% file: define5arch.tex
\documentclass[12pt]{amsart} 
\usepackage[latin1]{inputenc}
\usepackage[T1]{fontenc}

\usepackage{graphicx}
\usepackage{url}
\usepackage{inputenc} 
\usepackage{latexsym}
\usepackage{amssymb}
\usepackage{amsmath}
\usepackage{amsfonts}

\makeatletter
\def\inputonly#1{\inputonly@i#1,,\@nil}
\def\inputonly@i#1,#2,#3\@nil{%
  \global\@namedef{#1input}{#1}
  \ifx\relax#2\relax\else\inputonly@i#2,#3\@nil\fi}
\def\Input#1{\expandafter\ifx\csname#1input\endcsname\relax\else\input{#1}\fi}
\makeatother
\inputonly{textver2}
\usepackage[all,poly]{xy}

\newtheorem{theorem}{Theorem}[section]
\newtheorem{proposition}[theorem]{Proposition}
\newtheorem{lemma}[theorem]{Lemma}
\newtheorem{corollary}[theorem]{Corollary}

\theoremstyle{definition}
\newtheorem{definition}[theorem]{Definition}

\theoremstyle{remark}

\newtheorem{remark}[theorem]{Remark}

\newcommand{\SL}{\mathrm{SL}}
\newcommand{\Ad}{\mathrm{Ad}}
\newcommand{\GL}{\mathrm{GL}}
\newcommand{\PGL}{\mathrm{PGL}}

\begin{document}

\title[Convex projective polyhedral reflection groups ]{The definability criterions for convex projective polyhedral reflection groups }







\author{Kanghyun Choi}
\address{Sang-dong, Sarang Maeul Apt. 1612-1001, 
Buheung-ro, Wonmi-gu, Bucheon-si, Gyeonggi-do, Republic of Korea }
\email{kchoi1982@kaist.ac.kr}

\author{Suhyoung Choi }
\address{Department of Mathematical Sciences, KAIST, Daejeon 305-701, Republic of Korea}
\email{schoi@math.kaist.ac.kr}
\date{\today}

\begin{abstract}
Following Vinberg, we find the criterions for a subgroup generated by reflections 
$\Gamma \subset \SL^{\pm}(n+1,\mathbb{R})$ 
and its finite-index subgroups 
to be  definable over $\mathbb{A}$ where $\mathbb{A}$ is an integrally closed Noetherian ring in the field $\mathbb{R}$. 
We apply the criterions for  groups generated by reflections that act cocompactly on  
irreducible properly convex open subdomains of the $n$-dimensional projective sphere. 
This gives a method for constructing  injective group homomorphisms from such  Coxeter groups to
$\SL^{\pm}(n+1,\mathbb{Z})$. Finally we provide some examples of 
$\SL^{\pm}(n+1,\mathbb{Z})$-representations of such Coxeter groups. 
In particular, we consider  simplicial reflection groups that are isomorphic to hyperbolic simplicial groups 
and classify all the conjugacy classes of 
the reflection subgroups in $\SL^{\pm}(n+1,\mathbb{R})$ that are definable over $\mathbb{Z}$.
These were known by Goldman, Benoist, and so on previously. 
\end{abstract}


\thanks{Both authors were 
supported by Mid-career Researcher Program through NRF grant (2011-0027600) funded by the MEST}

\maketitle

\tableofcontents


\input{textver2}

\end{document}

%% file: textver2.tex
\section{Introduction}
\label{sec:intro}

Let $\mathbb{F}$ be a field of characteristic $0$ and let $V$ be an (n+1)-dimensional vector space over $\mathbb{F}$. We define
  \[\SL^{\pm}(n+1,\mathbb{F}):=\{A \in \GL(n+1,\mathbb{F}) | \det A = \pm 1\}\]
  as the group of linear automorphisms of $V$ whose determinants are either $1$ or $-1$. 
  A {\em reflection} $\sigma$ is an element of order 2 of $\SL^{\pm}(n+ 1,\mathbb{F})$ which is the identity on a hyperplane.  
  Every reflection $\sigma$ acts on $V$ by
  \[\sigma_{\alpha, \vec{v}}:\vec{x} \rightarrow \vec{x}- \alpha(\vec{x})\vec{v}.\]
  Namely, all reflections are of the form $\sigma:=\sigma_{\alpha, \vec{v}}:=I -\alpha \otimes \vec{v}$ for some $\alpha \in V^{*}$ and $\vec{v}\in V$ with $\alpha(\vec{v})=2$.
  A reflection itself is determined up to an automorphism 
  \begin{equation} \label{eqn:amb}
  \vec{v}\rightarrow c\vec{v}, \alpha \rightarrow c^{-1} \alpha \hbox{ for } c \ne 0.
  \end{equation}   
  
  \raggedbottom
  
  
  Let $\Gamma$ be the group generated by reflections 
  $\sigma_{i}=\sigma_{\alpha_{i},\vec{v}_{i}}$, $i \in \mathbb{N}$. We set $c_{i,j}=\alpha_{i}(\vec{v}_{j})$ and 
  call the matrix $C=(c_{i,j})$ the \emph{Cartan matrix} of $\Gamma$.
  The Cartan matrix is determined by $\Gamma$ up to conjugations by positive diagonal matrices by the above ambiguity. 
  Products of the form $c_{i_{1}i_{2}}c_{i_{2}i_{3}}...c_{i_{k}i_{1}}$
   are called \emph{cyclic products}. A cyclic product is called $\emph{simple}$ 
   if all its indices $i_{1},...,i_{k}$ are distinct. It is clear that every cyclic product is a product of simple ones
   and are invariants of the Cartan matrix under the conjugations by the positive diagonal matrices. 

Let $O$ denote the origin in $\mathbb{R}^{n+1}$ from now on.
Let $\mathbb{S}^{n}$ be the $n$-dimensional real projective sphere, i.e. the quotient space $(\mathbb{R}^{n+1}-\{O\})/\sim$ where 
$\sim$ is the equivalence relation given by $\vec{v}\sim \vec{w}$ if and only if $\vec{v}=r\vec{w}$ for a positive real number $r$. 
The elements of $\SL^{\pm}(n+1, \mathbb{R})$ restrict to automorphisms of $\mathbb{S}^{n}$, which are the {\em projective automorphisms}. 
Then we can identify $\SL^{\pm}(n+1,\mathbb{R})$ with the group of projective automorphisms of $\mathbb{S}^{n}$.
Denote by $\Pi: \mathbb{R}^{n+1}-\{O\} \rightarrow {\mathbb{S}}^{n}$ the quotient map.
A {\em projective reflection} is the induced projective automorphism ${\mathbb{S}}^{n}$ induced by a reflection. 
(Sometimes, we will drop the adjective ``projective".)

A {\em cone} in $\mathbb{R}^{n+1}$ is a subset $C$ of $\mathbb{R}^{n+1}-\{O\}$ so that $v \in C$ iff $s v \in C$ for $s > 0$. 
A {\em polyhedral cone} is a cone with nonempty interior bounded by finitely many hyperplanes. 
Let $P \subset \mathbb{S}^{n}$ be an $n$-\emph{dimensional convex polytope}, the image under $\Pi(C)$ in  $\mathbb{S}^{n}$ 
of a convex closed polyhedral cone $C$ of $\mathbb{R}^{n+1}-\{O\}$ with nonempty interior. 

Vinberg ~\cite{Vinberg2} initiated a ground-breaking study to construct injective homomorphisms from a Coxeter group 
that  is based on an $n$-dimensional polytope to $\SL^{\pm}(n+1,\mathbb{R})$: 
the group generated by projective reflections by sides of an $n$-dimensional polytope in 
$\mathbb{S}^{n}$. The properties of the corresponding representation spaces are recently studied by 
Benoist(\cite{Benoist1}, \cite{Benoist2}, \cite{Benoist3}, \cite{Benoist4}), Choi(\cite{Choi2}), and Marquis(\cite{Marquis}). 
See Proposition 6.2 and Figures 1 -- 4 in Section 6 for the parametrized subspaces in the case of triangular reflection groups.

For projective geometry, one normally works with the projective linear group  $\PGL(n+1, \mathbb{R})$ 
acting on the real projective space ${\mathbb{R}}P^n$.
   We work with $\SL^{\pm}(n+1,\mathbb{R})$ as this is a more standard Lie group to study the arithmetic questions.
   (See some discussion at Chapter 3 of \cite{psconv}, Section 2 of \cite{Cooper1} and \cite{Cooper}.)

A {\em great circle} in $\mathbb{S}^{n}$ is the image $\Pi(V -\{O\})$ of a two-dimensional subspace $V$ in $\mathbb{R}^{n+1}$.
A {\em segment} is a proper connected subset of a great circle.  A {\em convex} segment is one that does not contain an antipodal pair of
points except at the end points. A {\em closed hemisphere} is the image $\Pi(H - \{O\})$ where $H$ is a subset of $\mathbb{R}^{n+1}$
given by a linear inequality $f(\vec v) \geq 0$. A {\em hypersphere} is the image of $\Pi(W-\{O\})$ for a subspace $W$ of $\mathbb{R}^{n+1}$.
An {\em open hemisphere} is the interior of a closed hemisphere.
A {\em convex subset} of $\mathbb{S}^{n}$ is a set such that any pair of its point is connected by a convex segment.
A convex subset is either $\mathbb{S}^{n}$ itself or a subset of a closed hemisphere given by a linear inequality. 
(See Proposition 2.3 of \cite{psconv}.) 
A {\em properly convex subset} is a convex subset whose closure does not contain a pair of antipodal points.  
It is a precompact subset of an open hemisphere $H$. (See Proposition 2.4 of \cite{psconv}.) 
We will require a convex polytope to be properly convex always in this paper. 
A properly convex domain is {\em strictly convex} if the boundary contains no nontrivial segment.

 Let $\Omega$ be a properly convex open domain in
 $\mathbb{S}^{n}$ and $\Gamma$ is a discrete subgroup of $\SL^{\pm}(n+1,\mathbb{R})$ acting on $\Omega$ properly discontinuously.  
 When the properly convex projective orbifold $Q=\Omega/\Gamma$ is compact and Hausdorff, 
 we say that $\Gamma$ \emph{divides} $\Omega$. 
A group $\Gamma \subset \SL^{\pm}(n+1,\mathbb{R})$ generated by projective reflections 
that divides a properly convex open domain $\Omega$ in $\mathbb{S}^{n}$ 
is called a {\em dividing projective  reflection group}.
In this case, there exists a (properly) convex polytope $P$ so that $\Gamma$ is generated by 
projective reflections fixing sides of $P$. Hence, $\Gamma$ is called a {\em dividing polyhedral projective reflection group}.

Let $O^+(1, n)$ denote the subgroup of  $\SL^{\pm}(n+1,\mathbb{R})$ preserving the quadratic form 
\[q(x_0, x_1, \dots, x_n) := x_0^2 - x_1^2 - \cdots - x_n^2\] on $\mathbb{R}^{n+1}$ with coordinates
$x_0, x_1, \dots, x_n$ and preserving the positive cone of $q> 0$. 
The positive standard cone $C$ given by $q>0, x_0 > 0$ is projected to the interior $B$ of 
a conic. $B$ is called the {\em standard ball} in $\mathbb{S}^n$. This forms a Klein model of 
a hyperbolic $n$-space as $B$ admits a complete hyperbolic metric that is invariant under the action of $O^+(1, n)$. 
For any a cocompact lattice $\Gamma$ in $O^+(1, n)$, 
$B/\Gamma$ is compact Hausdorff, and 
$\Gamma$ is an example of a dividing polyhedral projective reflection group. 

There are numerous examples of dividing polyhedral projective reflection groups that are not in $O^+(1, n)$
obtained by Benoist \cite{Benoist4}, Choi \cite{Choi}, and Marquis \cite{Marquis}. 
(See also \cite{CHL} and \cite{CL}.)

Let $\Delta$ be a family of linear automorphisms of $V$.
A ring   $\mathbb{A}$ in $\mathbb{F}$ is called a \emph{ring of definition} for $\Delta$ if $V$ contains an $\mathbb{A}$-lattice
that is invariant under $\Delta$. In that case we also say that $\Delta$ is \emph{definable over} $\mathbb{A}$.
If $\mathbb{A}$ is a principal ideal domain, then $\Delta$ is definable over $\mathbb{A}$ if and only if 
there is a basis in which the automorphisms of $\Delta$ can be written down by matrices with elements in $\mathbb{A}$.
Our first main result is Theorem \ref{thm:mainprop} following from Theorem \ref{thm:main2} in Section 5. 
This generalizes Theorem 5 of Vinberg \cite{Vinberg3} for the case of subgroups of $O^+(1, n)$. (See also \cite{Vinberg1}.)

\begin{theorem}\label{thm:mainprop}
Suppose that $n\geq 2$ is an integer. Let $\Gamma \subset \SL^{\pm}(n+1,\mathbb{R})$ be 
a polyhedral projective  reflection group
dividing a strictly convex open domain $\Omega$. 
Let $\Delta$ be any subgroup of finite index in $\Gamma$.
Then the following 
properties of an integrally closed Noetherian ring $\mathbb{A}\subset \mathbb{R}$ are equivalent.
\begin{enumerate}
\item[{\rm (a)}] $\Gamma$ is definable over $\mathbb{A}$.
\item[{\rm (a${}'$)}] $\mathbb{A}$ contains all the simple cyclic products of the Cartan matrix of $\Gamma$.
\item[{\rm (b)}] $\Delta$ is definable over $\mathbb{A}$.
\end{enumerate}

\end{theorem}
The equivalence of (a) and (a${}'$) is shown by Vinberg \cite{Vinberg3}, thanks to 
the Zariski density theorem of Benoist \cite{Benoist0}. 
Kac and Vinberg were first to find examples of triangular projective reflection groups in  $\SL^{\pm}(3,\mathbb{Z})$  (~\cite{KV}).
Our main result is their equivalence to (b).

See Section \ref{sec:examples} and Proposition \ref{prop:triangle} for concrete examples of such group actions.
For ``irreducible domain'' see Subsection \ref{sec:reporb} or Benoist \cite{Benoist2}. 
A strictly convex open domain is an irreducible domain.
Let $\Gamma$ be as in the premise of the above theorem. 
According to Vey \cite{Vey}, 
$\Gamma$ is necessarily strongly irreducible, i.e., 
a finite index subgroup of the linear group $\Gamma$ are irreducible since $\Omega$ is strictly convex and hence is irreducible.
(See Theorem 1.1 of \cite{Benoist2}.)
Conversely, if $\Gamma$ is not a product group virtually or if $\Gamma$ is strongly irreducible,
then $\Omega$ has to be irreducible also.

Our strategy to prove Theorem \ref{thm:mainprop} is  to prove that the theorem holds for $\mathbb{F}$-Zariski-dense projective reflection groups  
in $\SL^{\pm}(n+ 1,\mathbb{F})$
where $\mathbb{F}$ is an algebraically closed field of characteristic $0$; 
i.e., Theorem \ref{thm:main2}.
In the proof, we use the Vinberg theory on rings of definitions of subgroups of linear algebraic groups.
In particular, we made a use of the fact that if $\Delta$ is definable over $\mathbb{A}$ then  $\Ad(\Gamma)$ 
is definable over $\mathbb{A}$ .  
Along the theory of Vinberg, Lemma \ref{grinberg} is essential for the proof.
Theorem \ref{thm:main2}, the result of Vinberg for orthogonal groups, and the Zariski density theorem of Benoist (Theorem \ref{thm:zariski}) imply Theorem \ref{thm:mainprop}.

The equivalence of (a) and (b) in Theorem \ref{thm:mainprop} can be proved for subgroups generated by finitely many elements 
if we add mild conditions on generators.
These results are Theorem \ref{thm:maingeneral} 
which we state here.

\begin{theorem}\label{thm:main} 
Suppose $n\geq 2$ is an  integer, and let $\mathbb{F}$ be an algebraically closed field of characteristic $0$. 
Let $\Gamma \subset \SL^{\pm}(n+1,\mathbb{F})$ be an $\mathbb{F}$-Zariski-dense 
subgroup generated by finitely many elements $\{\sigma_1, \dots, \sigma_m\}$. 
Let $\Delta$ be any subgroup of finite index in $\Gamma$.  Let $\mathbb{A}$ be an integrally closed Noetherian ring in $\mathbb{F}$, 
and let $\mathbb{K}$ be the field of fractions of $\mathbb{A}$. Suppose that for each $\sigma_i$ in $\{\sigma_1, \dots, \sigma_m\}$ 
satisfies either one of the two following conditions{\rm :}
\begin{itemize}
\item $\mathrm{tr} \sigma_{i} \ne 0$ and $\mathrm{tr} \sigma_i \in \mathbb{K}$ or 
\item The order of $\sigma_{i}$ is finite and relatively prime to $n+1$. 
\end{itemize}
Then the following are equivalent.
\begin{enumerate}
\item[{\rm (a)}] $\Gamma$ is definable over $\mathbb{A}$.
\item[{\rm (b)}] $\Delta$ is definable over $\mathbb{A}$.

\end{enumerate}
\end{theorem}
Note that different $\sigma_i$s can satisfy different conditions. 
We are unable to determine yet if we can drop the conditions entirely.

Section 2 reviews basic definitions of real projective orbifolds and states the Zariski density theorem of Benoist, i.e., Theorem 1.3 of \cite{Benoist0}. 
Section 3 is a review of the the Vinberg theory of polyhedral projective reflection groups. Section 4 is that of the theory of
rings of definition of Zariski-dense subgroups of semisimple linear algebraic groups. Section 5 states and proves our main results.
Finally, in Section 6 we apply our main results of the paper 
to some dividing polyhedral reflection groups to determine if they are definable over $\mathbb{Z}$. 
In particular, we will consider dividing projective simplicial reflection groups that are isomorphic to hyperbolic simplicial groups and  classify 
all the conjugacy classes of projective reflection groups that are definable over $\mathbb{Z}$, an unpublished work of Goldman.
(See Proposition \ref{prop:triangle}, \ref{prop:tetrahedron}, and \ref{prop:4-simplex}.)
Our examples of dividing polyhedral projective reflection groups that are definable over $\mathbb{Z}$ include some triangular, tetrahedral, and cubical reflection groups. 
 Our results are also applicable to many hyperbolic Coxeter orbifolds including orderable ones \cite{Choi2} 
 (see \cite{Marquis}).

  We thank William Goldman for introducing us to this theory. 
  In fact, he has found some of these examples himself.
  We thank Alan Reid and Dave Witte Morris for giving us many informations on arithmetic groups,
 and Gye-Seon Lee for helpful discussions and for letting us use his graphics of a cube and a prism.
 We also thank Ja Kyung Koo and his students for having joint seminars on arithmetic hyperbolic geometry
 which helped us much.


\section{Real projective orbifolds}
\subsection{Orbifolds}
Given two manifolds $M_1$ and $M_2$ with groups $G_1$ and $G_2$ acting on them respectively,
a map $f:M_1 \rightarrow M_2$ is \emph{equivariant} with respect to a homomorphism $h:G_1 \rightarrow G_2$ if
$h(g) \circ f = f \circ g$ holds for each $g \in G_1$.

We will be using the language of orbifolds for later purposes. An $n$-dimensional \emph{orbifold structure} on
a second countable Hausdorff space $X$ is
given by atlas of compatible charts $(U, G, \phi)$ where $U$ is an open subset of $\mathbb{R}^n$, $G$ is a finite group acting on $U$,
and $\phi:U \rightarrow V$ induces a homeomorphism $U/G \rightarrow V$ for an open subset $V$ of $X$.
We say that $\phi(U)$ is \emph{modeled} on $(U, G)$ or $(U, G, \phi)$ also.
Two charts $(U, G, \phi)$ and $(U', G', \phi')$ are
\emph{compatible} if given any point $p \in \phi(U) \cap \phi(U')$, there is a chart $(U'', G'', \phi'')$ so that
$\phi''(U'')$ is an open neighborhood of $p$ in $\phi(U) \cap \phi(U')$ where there are smooth lifting embeddings
$\tilde \phi''_1:U'' \rightarrow U$ 
equivariant with respect to an injective homomorphism $G'' \rightarrow G$ and $\tilde \phi''_2: U'' \rightarrow U'$
equivariant with respect to one $G'' \rightarrow G'$.

Here, $X$ with an orbifold structure is said to be an \emph{orbifold} and $X$ is said to be the underlying space of the orbifold.

Given a manifold $M$ and a discrete group $\Gamma$ acting on $M$
properly discontinuously and but not necessarily freely, we can form $M/\Gamma$ as a quotient space.
$M/\Gamma$ has a natural orbifold structure given by covering $M$ by open sets $U$ so that
if $g(U) \cap U \ne \emptyset$, then $g(U) = U$. The orbifold structure is given by the collection of every
$(U, G, \phi)$ where $G$ is the finite subgroup of $\Gamma$ acting on a Euclidean open set $U$ of $M$. 

Since the orbifolds we study here admit real projective structures,
they are of the form $M/\Gamma$ for a manifold $M$ and a discrete group $\Gamma$ acting on $M$
properly discontinuously and but not necessarily freely by Thurston \cite{thnote}
as these are ``good orbifolds". (See Theorem 6.1.1 of \cite{Choi3}.) 
Two such orbifolds $M/\Gamma$ and $N/\Delta$ for simply connected manifolds $M$ and $N$ with discrete groups 
$\Gamma$ and $\Delta$ acting properly discontinuously
 are \emph{diffeomorphic} if
there is a diffeomorphism $f:M \rightarrow N$ equivariant with respect to an isomorphism $h:\Gamma \rightarrow \Delta$.
The \emph{fundamental group} $\pi_{1}(Q)$ is  the abstract group which is isomorphic to $\Gamma$
whenever $M$ is simply connected.  
If $M$ is simply connected, $M$ is said to be a \emph{universal cover} of $M/\Gamma$.
(For the general definition of orbifolds and geometric structures on it, we refer to \cite{Choi} and \cite{Choi3}.)

\subsection{Real projective orbifolds} \label{sec:reporb}
 Let $n\geq 2$ and $V$ be an $n$-dimensional real vector space $\mathbb{R}^{n+1}$. We identify
\begin{displaymath}
\SL^{\pm}(n+1,\mathbb{R}):=\{A \in \GL(n+1,\mathbb{R}) | \det A = \pm 1\}
\end{displaymath}
 with the group of projective automorphisms of $\mathbb{S}^{n}$.
 Recall that a subdomain $\Omega$ of  $\mathbb{S}^{n}$  is \emph{properly convex} 
 if it is convex and its closure $\overline{\Omega}$ does not contain two antipodal points.
 $\Omega$ is \emph{strictly convex} if furthermore its boundary $\partial \Omega$ does not contain any line segment of positive length.
 A strictly convex domain is properly convex. 

 A \emph{properly convex projective orbifold}, $Q$, is of the form $\Omega/\Gamma$ where $\Omega$ is a properly convex open domain in
 $\mathbb{S}^{n}$ and $\Gamma$ is a discrete subgroup of $\SL^{\pm}(n+1,\mathbb{R})$  acting on $\Omega$ properly discontinuously.
 A \emph{properly convex projective structure} on an orbifold $M$ is a diffeomorphism $f: M \rightarrow Q$ for $Q$ as above.

 The two properly convex projective orbifolds
 $Q_{1}=\Omega_{1}/\Gamma_{1}$ and $Q_{2}=\Omega_{2}/\Gamma_{2}$ are \emph{projectively diffeomorphic} if there exists a projective automorphism
 $h \in \SL^{\pm}(n+1,\mathbb{R})$ such that $h(\Omega_{1})= \Omega_{2}$ and $h \Gamma_{1}h^{-1}=\Gamma_{2}$. 

 Given a properly convex projective structure on an orbifold $M$, we obtain a homomorphism
 $h: \pi_1(M) \rightarrow \SL^{\pm}(n+1,\mathbb{R})$, called a {\em holonomy homomorphism} determined up to
 conjugation by $\SL^{\pm}(n+1,\mathbb{R})$: Let $\Omega/\Gamma$ be a properly convex projective orbifold
 and $f: M \rightarrow \Omega/\Gamma$ be a projective diffeomorphism. Then let $f_*: \pi_1(M) \rightarrow \Gamma$
 be the induced homomorphism. Also $f_*$ considered as a map  $\pi_1(M) \rightarrow \SL^{\pm}(n+1,\mathbb{R})$
 is the {\em holonomy homomorphism}. The image of the holonomy homomorphism is called the {\em holonomy group} of $M$ or $\pi_1(M)$.
 
 The usual definition of a properly convex real projective orbifold as in \cite{JM} or \cite{Choi2} is equivalent to this. (See Section 3.4 of \cite{Thbook} 
 or Section 2.2 of \cite{Choi2} for discussions.) 

 Any open domain $\Omega$ in $\mathbb{S}^{n}$ is the image of the unique convex open cone $C$ in $V$ under the projection map.
 The open convex cone $C$ in $V-\{O\}$ is said to be \emph{reducible} if it can be written as the {\em sum} 
 \[C = C_{1} \oplus C_{2}: = \{ v+w | v \in C_1, w \in C_2\}\]
 of two convex cones $C_{i}$ in proper subspaces $V_{i}$ of $V$ where 
 $V$ is a direct sum $V_1 \oplus V_2 $. 
 A properly convex open domain $\Omega$ in $\mathbb{S}^{n}$ is said to be \textit{reducible}
 if its preimage $C$ is reducible. Otherwise we say that $\Omega$ is \emph{irreducible}.
 A strictly convex domain is irreducible since a reducible domain always have a nontrivial segment in the boundary. 
 
 Let $x_0, x_1, \dots, x_{n}$ denote the set of standard coordinates of the vector space $\mathbb{R}^{n+1}$.
 The \emph{hyperbolic $n$-space}  $\mathbb{H}^{n} = B \subset \mathbb{S}^n$ is defined as the image $\Pi(\Lambda_{n+1})$ of 
 the positive Lorentz cone
 \[\Lambda_{n+1} :=\{x \in \mathbb{R}^{n+1}| q(x) > 0 \mbox{ and } x_{1} > 0 \}\] 
 where $q(x) = x_{0}^{2}-x_{1}^{2}- \cdots -x_{n}^{2}$ holds. 
 The group of projective automorphisms of $\mathbb{H}^n$ is the group
 \[\mathrm{Aut}(\mathbb{H}^{n}) = \mathrm{O}^{+}(1,n) \subset \SL^{\pm}(n+1,\mathbb{R})\] 
 of orthogonal transformations of $q$ which preserve $\mathbb{H}^{n}$.
 If there exists a discrete subgroup $\Gamma$ of $\mathrm{O}^{+}(1,n)$ which preserves $\mathbb{H}^{n}$, 
 we call the quotient projective orbifold
 $\mathbb{H}^{n}/\Gamma$ a $\emph{hyperbolic orbifold}$. 
 Note that $\mathbb{H}^{n}$ is an example of  an irreducible strictly convex open domain in
 $\mathbb{S}^{n}$. 
 

\subsection{The Zariski density theorem of Benoist}

Now we state our version  of the theorem which follows directly from the version of Cooper and Delp. (Corollary 4.2 of \cite{Cooper})
\begin{theorem}[Theorem 1.3 \cite{Benoist0}]\label{thm:zariski}
Let $\Gamma$ be a discrete subgroup of $\SL^{\pm}(n+1,\mathbb{R})$ which divides an irreducible strictly convex open domain
$\Omega$ in $\mathbb{S}^{n}$. The Zariski closure $\overline{\Gamma}$ over $\mathbb{R}$ is either $\mathrm{O}^{+}(1,n)$ 
or $\mathrm{SO}^{+}(1,n)$ iff $\Omega$ is $\mathbb{H}^{n}$.
Otherwise $\overline{\Gamma}$ is either $\SL^{\pm}(n+1,\mathbb{R})$ or $\SL(n+1,\mathbb{R})$.
\end{theorem}
\begin{proof}
 Let $P \colon \mathbb{S}^{n}  \to \mathbb{R}P^{n}$ denote the double covering map. 
 If a discrete subgroup $\Gamma$ of $\SL^{\pm}(n+1,\mathbb{R})$ divides an irreducible strictly convex open domain
$\Omega$ in $\mathbb{S}^{n}$, then $\Gamma$ divides $P(\Omega)$ in $\mathbb{R}P^{n}$. 
Now the result follows from Corollary 4.2 of \cite{Cooper}. 
\end{proof}

A group $\Gamma$ of $\SL^{\pm}(n+1,\mathbb{R})$ is said to be \emph{irreducible}
 if there are no $\Gamma$-invariant nontrivial subspaces in $\mathbb{R}^{n+1}$, 
 and is said to be \emph{strongly irreducible} if the finite-index subgroups are irreducible.

\begin{definition}
A discrete subgroup of $\SL^{\pm}(n+1,\mathbb{R})$ which divides an irreducible properly convex open domain
$\Omega$ in $\mathbb{S}^{n}$ is said to be an {\em irreducible dividing group}.
\end{definition}

Of course this group is a strongly irreducible linear subgroup of $\SL^{\pm}(n+1,\mathbb{R})$
by the Vey irreducibility theorem (Theorem 5.1 of ~\cite{Benoistsurvey}).

\section{Polyhedral projective reflection groups}

\subsection{The Vinberg Condition}
A {\em rotation} $g$ is an element of $\SL^{\pm}(n+1,\mathbb{R})$ that fixes all points of 
a codimension-two subspace $W$ in ${\mathbb{R}}^{n+1}$ and 
can be written an orthogonal transformation of the form
\[
\left(
\begin{array}{cc}
\cos \theta  & \sin \theta    \\
-\sin \theta  & \cos \theta  
\end{array}
\right), 0 < \theta \leq \pi
\]
with respect to a basis in the complement $W'$. 
We say that $g$ is a {\em rotation} about the subspace $\Pi(W - \{O\})$ of {\em angle} $\theta$.

A $k$-\emph{face} of a convex polytope $P$ is a $k$-dimensional convex subset
of $P$ obtained as an intersection of $P$ with some hyperspheres which do not meet the interior $P^{o}$.
A \emph{face} is an $(n-1)$-face. 
A \emph{projective reflection} is a reflection defined in $n$-dimensional real vector space, i.e. a reflection in
$\SL^{\pm}(n+1,\mathbb{R})$. Let $S$ be the set of faces of $P$ and for every $s$ in $S$, we can associate a projective reflection
$\sigma _{s} = I - \alpha _{s} \otimes \vec{v}_{s}$ with $\alpha_{s}(\vec{v}_{s}) = 2$ which fixes s.
A suitable choice of signs allows us to suppose that $P$ is defined by the inequalities 
\[\alpha_{s}\leq0 \hbox{ for } {s\in S}.\]

Let $\Gamma$ be the group generated by the reflections $\{\sigma_{s}|s \in S\}$. 
The following theorem of Vinberg ~\cite{Vinberg2} provides a necessary and sufficient condition for such a group $\Gamma$
to act on a convex subdomain of $\mathbb{S}^{n}$. Let $c_{s,t} := \alpha_{s}(\vec{v}_{t})$ for $s,t \in S$. 
\begin{itemize}
\item $c_{s,t} c_{t,s}=4\cos^{2}(\frac{\pi}{m_{s,t} })$ for an integer $m_{s, t} \geq 3$ is the condition that $\{\sigma_s, \sigma_t\}$ generate 
a  dihedral group of order $2m_{s,t}$ acting on $\mathbb{S}^n$ discretely fixing a codimension-two subspace
and 
\item $c_{s, t}=0, c_{t, s}=0$ are the condition that $\{\sigma_s, \sigma_t\}$ generate a dihedral group of
order $4$ acting on ${\mathbb{S}}^n$ fixing all points of a codimension-two subspace similarly. 
\end{itemize} 

Two faces of $P$ are {\em adjacent} if they meet at a codimension-two side in $P$.
A codimension-two side $u$ have the {\em edge order} or {\em order} assigned to be $m_{s, t}$ if the faces $s$ and $t$ meet at $u$

\begin{theorem}[Vinberg]\label{reflection}
Let $P$ be a convex compact polytope of $\mathbb{S}^{n}$ and, for each face $s$ of $P$,
let $\sigma _{s} = I - \alpha _{s} \otimes \vec{v}_{s}$ be a projective reflection fixing this face $s$.
Let $\Gamma$ be the group generated by the reflections $\sigma_{s}$.
Then the following conditions for every pair of faces $(s, t)$ are necessary and sufficient condition 
for $\Gamma$ to act on some convex subdomain $\Omega$ of $\mathbb{S}^{n}$ with the  fundamental domain $P$
\begin{enumerate}
\item[{\rm (L1)}] $c_{s,t} \leq 0$ 
\item[]
\begin{itemize}
\item[{\rm (L2(i))}] $c_{s,t}c_{t,s} \geq 4$ if $s$ and $t$ not adjacent; or 
\item[{\rm (L2(ii))}] $c_{s,t,}c_{t,s}=4\cos^{2}(\frac{\pi}{m_{s,t} })$ with  an integer $m_{s,t} \geq 2$
or \[c_{s,t}=0 \Leftrightarrow c_{t,s}=0\] if $s$ and $t$ are adjacent.
\end{itemize}
\end{enumerate}
Moreover $\Gamma$ is discrete. The convex domain $\Omega$ is open if and only if for every $x$ in $P$,
the group $\Gamma_{x}$ generated by $\sigma _{s}$ for $s$ containing $x$ is a finite group.
In this case, $\Gamma$ acts on $\Omega$  properly discontinuously with the compact quotient.
\end{theorem}
\begin{proof}
Proposition 17 of ~\cite{Vinberg2} gives the necessity of the condition (L1) and (L2).
\newline
Given (L1) and (L2), Theorem 1 and Proposition 6 and of  ~\cite{Vinberg2} 
show that $\Gamma$ preserves some convex subdomain $\Omega$ of $\mathbb{S}^{n}$ with
the fundamental domain $P$. 
\end{proof}
We will call the group generated by projective reflections fixing sides 
of some $n$-dimensional convex polytope $P$ of $\mathbb{S}^{n}$
the \emph{polyhedral projective reflection group}. Any polyhedral projective reflection group has the presentation
 \begin{displaymath}
 \langle s_{i}| (s_{i}s_{j})^{n_{ij}}=1\rangle,\mbox{ }n_{ij} \in \mathbb{N},\mbox{ } i, j=1, \ldots ,n,
 \end{displaymath}
 where $n_{ij}$ is defined for a subset of $\{1, \ldots, n\}^2$ and
 $n_{ij}$ is symmetric in $i, j$ and $n_{ii} = 2$. 
 An abstract group which has the above presentation is called a {\em Coxeter group with $n$ generators}.

\subsection{Coxeter orbifolds}


  A \emph{reflection} in an open subset $U$ of $\mathbb{R}^n$ is a transformation $U \rightarrow U$ of order two fixing a hyperspace in it.
  An $n$-\emph{dimensional Coxeter orbifold structure} on an $n$-dimensional convex polytope
  $P$ is an orbifold structure on $P$ where 
  \begin{itemize}
\item  each point of the interior of each face has
  a chart modeled on an open subset in $\mathbb{R}^n$ with a reflection acting on it and  
  \item each point of the interior of each side of codimension $2$
  has a chart modeled on an open subset of $\mathbb{R}^n$ with a dihedral group generated by reflections.
  \end{itemize} 
  The convex polytope $P$ with the Coxeter orbifold structure is denoted by $\hat P$ and is said to be an $n$-dimensional Coxeter orbifold
  (more precisely a Coxeter orbifold of type III in \cite{Davis2}). 
  
  In our cases, $\hat P = M/\Gamma$ for a simply-connected manifold $M$ and $\Gamma$ is a discrete group acting
  properly discontinuously since $\hat P$ admits a real projective structure.
  (See \cite{Choi} for details.)

  If $M$ is a convex domain in $\mathbb{S}^n$ on which a polyhedral projective reflection group $\Gamma$ acts cocompactly and properly discontinuously,
  we say that $\Gamma$ is a \emph{convex projective reflection group}.
  In this case, the orbifold $M/\Gamma$ has a \emph{convex projective structure}.
  If $\Gamma$ acts on the convex domain $\mathbb{H}^{n}$,
  then $\Gamma$ is a  \emph{hyperbolic reflection group}.
  In the later case, $M/\Gamma$ has a {\em hyperbolic structure}.
  

\begin{proposition}\label{prop:isom}
\begin{itemize}
\item Given two $n$-dimensional compact Coxeter orbifolds $M_1$ and $M_2$,
$M_1$ is diffeomorphic to $M_2$ if and only if $\pi_1(M_1)$ is isomorphic to $\pi_1(M_2)$.
{\rm (Davis) }.
\item A discrete group $\Gamma \subset \SL^\pm(n+1, \mathbb{R})$ generated 
by projective reflections divides a properly convex domain $\Omega$ in $\mathbb{S}^n$
if and only if it is a holonomy homomorphism of a compact properly convex projective Coxeter orbifold $M_1$.
\item Let $M_1$ be a compact Coxeter orbifold. 
The set of holonomy homomorphisms of  properly convex real projective structures  $M_1$ is identical
with the set of homomorphisms of $\pi_1(M_1) \rightarrow \SL^\pm(n+1, \mathbb{R})$ dividing properly convex domains.
\item The one-to-one correspondence also exists if we take quotients of both sets by the conjugation action by $\SL^\pm(n, \mathbb{R})$.
\end{itemize}
\end{proposition}
\begin{proof}
Let $P_1$ and $P_2$ be the convex polytopes that are the respective underlying spaces of $M_1$ and $M_2$.
By Example 7.14 of \cite{Davis}, $\pi_1(M_1)$ and $\pi_1(M_2)$ are geometric reflection groups. 
The Coxeter systems of geometric reflections groups are of type PM${}^n$ in the terminology of Davis by Corollary 8.2.10 of \cite{Davis}. 
(See Section 13.3 of \cite{Davis} for the definition of the Coxeter systems of type PM${}^n$. In brief, this means that the Coxeter system 
has the nerve of a pseudo-manifold. A geometric reflection group based on a convex polytope is clearly so since $\partial P$ is a topological manifold.)
Then the Coxeter systems give complete combinatorial data of the faces of $P_1$ and $P_2$
and the intersection pattern. 
If $\pi_1(M_1)$ is isomorphic to $\pi_1(M_2)$ and these have Coxeter systems 
of type PM${}^n$, then the combinatorial data of the faces of $P_1$ and $P_2$
are identical by Theorem 13.4.1 of \cite{Davis}.
Hence, there exists  a homeomorphism $P_1 \rightarrow P_2$ preserving the combinatorics. 
By Corollary 1.3 of \cite{Davis2}, we obtain a diffeomorphism $M_1 \rightarrow M_2$. 

For the second item,
the compact orbifold $\Omega/\Gamma$ has the fundamental group isomorphic to $\Gamma$ isomorphic to $\pi_1(\Omega/\Gamma)$.
Let $k: \pi_1(\Omega/\Gamma) \rightarrow \Gamma$ be the identity map. 
Therefore, $k$ is a holonomy homomorphism of a properly convex projective structure on $\Omega/\Gamma$. 

Conversely, if $k$ is a holonomy homomorphism of such a structure on $M_1$, then $k(\pi_1(M_1))$ acts on
a properly convex domain $\Omega$ corresponding to the universal cover of $M_1$ so that $\Omega/\Gamma$ is projectively
diffeomorphic to $M_1$. 

The remaining two items are consequences of the first two items.  
\end{proof}

\section{The results of Vinberg on rings of definitions}
In this section we state main results of Vinberg ~\cite{Vinberg3} regarding rings of definition of
Zariski-dense subgroups of semisimple linear algebraic groups. Our main result in the next section is heavily based on the Vinberg theory.

Let $V$ be a $(n+1)$-dimensional vector space over field $\mathbb{F}$ and $\mathbb{A}$ be a subring of $\mathbb{F}$.
A set $L$ in $V$ is called an $\mathbb{A}$-$lattice$ if it is a finitely generated $\mathbb{A}$-submodule and the natural map
$\mathbb{F}\bigotimes_{\mathbb{A}}L\rightarrow V$ is an isomorphism. If $\mathbb{A}$ is a principal ideal domain, then every
$\mathbb{A}$-lattice has a basis which at the same time is a basis of $V$ over $\mathbb{F}$.
If $\mathbb{K}$ is a subfield of $\mathbb{F}$, then a $\mathbb{K}$-lattice is just an $(n+1)$-dimensional vector space over $\mathbb{K}$.

Let $\Delta$ be a family of linear transformations of $V$. An integrally closed Noetherian ring 
$\mathbb{A}$ in $\mathbb{F}$ is called a \emph{ring of definition}
for $\Delta$ if $V$ contains an $\mathbb{A}$-lattice
that is invariant under $\Delta$. In that case we also say that $\Delta$ is \emph{definable over} $\mathbb{A}$. If $\mathbb{A}$ is a principal ideal domain,
then the fact that $\Delta$ is definable over $\mathbb{A}$ means
that there is a basis in which the transformations of $\Delta$ can be written down by matrices with elements in $\mathbb{A}$.
If $\Delta$ is definable over $\mathbb{A}$ and $\mathbb{B}$ is a ring that contains $\mathbb{A}$, then $\Delta$ is also definable over $\mathbb{B}$.

Now we assume that $\mathbb{F}$ is an algebraically closed field of characteristic $0$, $G$ a semisimple algebraic group over $\mathbb{F}$ 
in $\GL(n+1, \mathbb{F})$
and $\Gamma$ a subgroup Zariski-dense in $G$ over $\mathbb{F}$. Let $\Ad$ be the adjoint representation of $G$. Then $\Ad$ is a map from
$G$ to the automorphism group of its Lie algebra $\mathfrak{g}$. For any matrix group $H$, let $\mathrm{tr} H$ be the set of traces 
$\{\mathrm{tr} h | h \in H \}$. 

\begin{proposition}[Theorem 1 of ~\cite{Vinberg3}]\label{prop:thm1}
An integrally closed Noetherian ring $\mathbb{A}$ in $\mathbb{F}$ is a ring of definition for $\Ad (\Gamma)$ if and only if 
$\mathbb{A} \supset \mathrm{tr} \Ad (g)$ for all $g \in \Gamma$.
\end{proposition}


\begin{corollary}[Corollary to Theorem 1 of ~\cite{Vinberg3}]\label{cor}
There exists a smallest field of definition for the group $\Ad (\Gamma)$. If this field is an algebraic number field, then there exists
a smallest ring of definition for $\Ad (\Gamma)$.
\end{corollary}
\begin{proposition}[Theorem 2 of ~\cite{Vinberg3}]\label{prop:thm2}
 Let $\mathbb{A}$ be an integrally closed Noetherian ring in $\mathbb{F}$ and $\mathbb{K}$ its field of fractions. 
 Then $\Gamma$ is definable over $\mathbb{A}$ if and only if it is definable over $\mathbb{K}$ and
 $\Ad (\Gamma)$ is definable over $\mathbb{A}$.
\end{proposition}

\begin{proposition}[Theorem 3 of ~\cite{Vinberg3}]\label{prop:thm3}
Let $\Gamma_{1}$ be a finite index subgroup of $\Gamma$. 
Then the classes of rings of definition for $\Ad (\Gamma_{1})$ and $\Ad (\Gamma)$ are identical.
\end{proposition}

Let $V$ be an $(n+1)$-dimensional vector space over a (not necessarily algebraically closed) field $\mathbb{F}$ and $V^{*}$ it dual. 
Let $\Gamma$ be a group generated by reflections. We denote by $\mathbb{Z}[\mathrm{tr} \Gamma]$ the ring with a unit element generated 
by the set $\mathrm{tr} \Gamma$ in $\mathbb{F}$.

\begin{proposition}[Lemma 11 of ~\cite{Vinberg3}]\label{prop:cyclic}
$\mathbb{Z}[\mathrm{tr} \Gamma]$ is the ring with a unit element generated by the simple cyclic products.
\end{proposition}
 A group $H$ in $\SL^{\pm}(n+ 1,\mathbb{F})$ is \emph{absolutely irreducible if} if it is irreducible as a
matrix group in $\SL^{\pm}(n+ 1,\bar{\mathbb{F}})$ where $\bar{\mathbb{F}}$ is the algebraic closure of $\mathbb{F}$.


\begin{proposition}[Lemma 12 of ~\cite{Vinberg3}]\label{prop:tr}
Suppose that the characteristic of the field $\mathbb{F}$ is $0$ and $\Gamma$ is absolutely irreducible.
An integrally closed Noetherian ring $\mathbb{A}$ in $\mathbb{F}$ is a ring of definition for $\Gamma$
if and only if $\mathbb{A}$ contains $\mathbb{Z}[\mathrm{tr} \Gamma]$.
\end{proposition}

Suppose that a nondegenerate scalar product $(,)$ is defined in $V$. 
The group of automorphisms of $V$ that preserve this scalar product $(,)$ 
is denoted by $\mathrm{O}(V)$ and its unimodular subgroup by $\mathrm{SO}(V)$.

\begin{proposition}[Theorem 5 of ~\cite{Vinberg3}]\label{thm5}
Suppose $n\geq 3$, and let $\mathbb{F}$ be an algebraically closed field of characteristic $0$. Let $\Gamma \subset \mathrm{O}(V)$ be a subgroup generated by reflections
Zariski-dense over $\mathbb{F}$.  
Let $\Delta$ be any subgroup of finite index in $\Gamma$. 
Then the following properties of an integrally closed Noetherian ring $\mathbb{A}\subset \mathbb{F}$ are equivalent.
\begin{enumerate}
\item[{\rm (a)}] $\Gamma$ is definable over $\mathbb{A}$.
\item[{\rm (a${}'$)}] $\mathbb{A}$ contains all the simple cyclic products of the Cartan matrix of $\Gamma$.
\item[{\rm (b)}] $\Delta$ is definable over $\mathbb{A}$.
\item[{\rm (c)}] $\Ad(\Delta)$ is definable over $\mathbb{A}$.
\item[{\rm (d)}] $\Ad(\Gamma)$ is definable over $\mathbb{A}$.
\end{enumerate}
\end{proposition}

\section{The definability results}

We will prove Theorem \ref{thm:mainprop} in this section.
Throughout this section, let $\mathbb{F}$
be an algebraically closed field of characteristic $0$. To prove Theorem \ref{thm:mainprop}, 
we prove that the similar result holds for finitely generated groups if we add some conditions.

\subsection{The definability under assumptions}

\begin{theorem}\label{thm:maingeneral}
Suppose $n\geq 2$ is an  integer, and let $\mathbb{F}$ be an algebraically closed field of characteristic $0$. 
Let $\Gamma \subset \SL^{\pm}(n+1,\mathbb{F})$ be a $\mathbb{F}$-Zariski-dense 
subgroup generated by finitely many elements $\sigma_{i}$ for $i=1, \cdot \cdot \cdot, m$. 
Let $\Delta$ be any subgroup of finite index in $\Gamma$.  Let $\mathbb{A}$ be an integrally closed Noetherian ring in $\mathbb{F}$, 
and let $\mathbb{K}$ be the field of fractions of $\mathbb{A}$. 
Suppose that for each $i=1, \cdot \cdot \cdot, m$, $\sigma_{i}$ satisfies one of the following conditions\,{\rm :}
\begin{itemize}
\item [{\rm (a)}] $\mathrm{tr} \sigma_{i} \in \mathbb{K}, \mathrm{tr} \sigma_i \ne 0$.
\item [{\rm (b)}]The order of $\sigma_{i}$ is finite and relatively prime to $n+1$.
\end{itemize}
Then the following properties are equivalent.
\begin{enumerate}
\item[{\rm (D1)}] $\Gamma$ is definable over $\mathbb{A}$.
\item[{\rm (D2)}] $\Delta$ is definable over $\mathbb{A}$.

\end{enumerate}
\end{theorem}

We consider the Lie algebra $\mathfrak{sl}_{n+1}\left(\mathbb{R}\right)$ as a subalgebra of the matrix algebra 
$\mathbf{M}_{n+1} \left(\mathbb{R}\right)$ which is realized from the decomposition  
\[\mathbf{M}_{n+1}\left(\mathbb{R}\right)=\mathfrak{sl}_{n+1}\left(\mathbb{R}\right) +\mathbb{R}\cdot I_{n+1}\]
 where $I_{n+1}$ is the $(n+1)\times(n+1)$ identity matrix.  
We need the following lemma of Grinberg \cite{Gr} with his proof. 
\begin{lemma}\label{grinberg}
 Let $\mathbb{L}$ be a field of characteristic $0$, and $\mathbb{K}$ be a subfield of $\mathbb{L}$. 
 Let $n\in\mathbb N$. Let $U\in \GL(n,\mathbb{L})$ be such that every $V\in\mathfrak{sl}_{n+1}\left(\mathbb{K}\right)$ 
 satisfies $UVU^{-1}\in \mathbf{M}_{n+1}\left(\mathbb{K}\right)$. Then, there exists a nonzero $\lambda\in \mathbb{L}$ such that $\lambda U \in \GL(n,\mathbb{K})$.
\end{lemma}
\begin{proof}
 
 Since $\mathbb{K}$ has characteristic $0$, we have $\mathbf{M}_{n+1}\left(\mathbb{K}\right)=\mathfrak{sl}_{n+1}\left(\mathbb{K}\right)+\mathbb{K}\cdot I_{n+1}$.
 We obtain
 \begin{enumerate}
\item[\rm (M)] every $V\in \mathbf{M}_{n+1}\left(\mathbb{K}\right)$ satisfies $UVU^{-1}\in \mathbf{M}_{n+1}\left(\mathbb{K}\right)$.
\end{enumerate}

Now, consider the map $r : \mathbf{M}_{n+1}\left(\mathbb{K}\right) \to \mathbf{M}_{n+1}\left(\mathbb{K}\right)$ which maps every $V$ to $UVU^{-1}$. 
This $r$ is well-defined due to (M), and is a $\mathbb{K}$-algebra isomorphism; hence, $r$ is a $\mathbb{K}$-algebra automorphism of 
$\mathbf{M}_{n+1}\left(\mathbb{K}\right)$. 
By the Skolem-Noether theorem, 
there exists some $P\in\GL(n,\mathbb{K})$ such that every $V\in \mathbf{M}_{n+1}\left(\mathbb{K}\right)$ satisfies $r\left(V\right) = PVP^{-1}$. 
So every $V\in\mathbf{M}_{n+1}\left(\mathbb{K}\right)$ satisfies
\[PVP^{-1} = r\left(V\right) = UVU^{-1}.\]
We rewrite this as $U^{-1}PV\left(U^{-1}P\right)^{-1}=V$. In other words, we have $U^{-1}PV=VU^{-1}P$. 
Since this holds for all $V\in \mathbf{M}_{n+1}\left(\mathbb{K}\right)$, it must also hold for all $V\in \mathbf{M}_{n+1}\left(\mathbb{L}\right)$ 
(because it is a linear equation in $V$, so it is enough to check it on an $\mathbb{L}$-basis of $\mathbf{M}_{n+1}\left(\mathbb{L}\right)$, 
but such a basis can be chosen to be in $\mathbf{M}_{n+1}\left(\mathbb{K}\right)$). In other words, the matrix $U^{-1}P$ is in the center of 
$\mathbf{M}_{n+1}\left(\mathbb{L}\right)$. However, this center is known to be $\mathbb{L}\cdot I_{n+1}$. Thus, we obtain $U^{-1}P\in \mathbb{L}\cdot I_{n+1}$. 
In other words, there exists some $\lambda\in \mathbb{L}$ such that $U^{-1}P = \lambda I_{n+1}$. 
This $\lambda$ is nonzero (else, $P$ would be $0$, 
contradicting $P\in \GL(n,\mathbb{K})$), so this becomes $P=\lambda U$. 
Hence, we obtain $\lambda U=P\in \GL(n,\mathbb{K})$.
(See also Theorem 4 of \cite{Vinberg3} for similar ideas.) 
\end{proof}


\begin{lemma}\label{lem:D2}
{\rm (D2)} implies that $\Gamma$ is definable over $\mathbb{K}$. 
\end{lemma}
\begin{proof} 
Assume (D2).   Then there exists an invertible linear map $g$ from $\mathbb{F}^{n+1}$ to $\mathbb{F}^{n+1}$ 
such that $\mathbb{K}^{n+1}$ is an invariant $\mathbb{K}$-form under $g \Delta g^{-1}$.    Then the set
\begin{displaymath}
\mathfrak{sl}_{n+1}\left(\mathbb{K}\right)=\{x \in \mathfrak{sl}_{n+1}\left(\mathbb{F}\right) | x \in \mathrm{End} \mathbb{K}^{n+1}\}
\end{displaymath}
is a $\mathbb{K}$-form of $\mathfrak{sl}_{n+1}\left(\mathbb{F}\right)$  invariant under $\Ad (g \Delta g^{-1})$.   
By the proof of Theorem 3 of  ~\cite{Vinberg3}, $\mathfrak{sl}_{n+1}\left(\mathbb{K}\right)$ is also invariant under $\Ad(g \Gamma g^{-1})$.                                                     
For a linear transformation $y$ from $\mathbb{F}^{n+1}$ to $\mathbb{F}^{n+1}$, let $\hat{y}$ denote the matrix of $y$ with respect to the standard basis of $\mathbb{F}^{n+1}$.  
Then since $\mathfrak{sl}_{n+1}\left(\mathbb{K}\right)$ is invariant under  $\Ad(g \Gamma g^{-1})$,  we obtain
\[\hat{g} \hat{\sigma_{i}} \hat{g}^{-1}x\hat{g} \hat{\sigma_{i}} \hat{g^{-1}} \in \mathbf{M}_{n+1}\left(\mathbb{K}\right)
\hbox{ for every } x \in \mathfrak{sl}_{n+1}\left(\mathbb{K}\right).\]  
By Lemma \ref{grinberg}, there exists a nonzero $\lambda_{i} \in \mathbb{F}$ such that $\lambda_{i} \hat{g} \hat{\sigma_{i}} \hat{g}^{-1}$ is in $\GL(n+1,\mathbb{K})$ so that 
\begin{equation}
\mathrm{det} (\lambda_{i} \hat{g}\hat{\sigma_{i}} \hat{g}^{-1})= \pm \lambda_{i}^{n+1} \in \mathbb{K} 
\tag{I}
\end{equation}
and
\begin{equation}
 \mathrm{tr} \lambda_{i} \hat{g} \hat{\sigma_{i}} \hat{g}^{-1}=\lambda_{i}\mathrm{tr} \sigma_{i} \in  \mathbb{K}.
\tag{II}
\end{equation}
 If $\sigma_{i}$ is of finite order $l_{i}$, we obtain  
 \[ (\lambda_{i}\hat{g}\hat{\sigma_{i}} \hat{g}^{-1})^{l_{i}} = \lambda_{i}^{l_{i}}I_{n+1} = I_{n+1} \hbox{ and hence } 
\lambda_{i}^{l_{i}} \in \mathbb{K}.
\tag{III}\]


Suppose that $\sigma_{i}$ satisfies condition (\rm a). From (II) and that $\mathrm{tr} \sigma_{i}$ is in  $\mathbb{K}$ and is invertible, 
we conclude that $\lambda_{i}$ is in $\mathbb{K}$. Since $\lambda_{i} \hat{g} \hat{\sigma_{i}} \hat{g}^{-1}$ is in 
$\GL(n+1,\mathbb{K})$, we conclude that 
\[\hat{g} \hat{\sigma_{i}} \hat{g}^{-1} \in \SL^{\pm}(n+1, \mathbb{K}).\]

Suppose that  $\sigma_{i}$ satisfies condition (\rm b) and let $l_{i}$ be the order of  $\sigma_{i}$ which is relatively prime to $n+1$.  
Then from (I) and (III) and that $l_{i}$ is relatively prime to $n+1$, we obtain 
\[\lambda_{i} \in \mathbb{K} \hbox{ and } \hat{g} \hat{\sigma_{i}} \hat{g}^{-1}  \in \SL^{\pm}(n+1, \mathbb{K}).\]
Therefore, $\Gamma$ is definable over $\mathbb{K}$. 
\end{proof}

\begin{proof}[The proof of Theorem \ref{thm:maingeneral}] 
The fact that $(D1)$ implies $(D2)$ is obvious.  
We consider the statements,
\begin{itemize}
\item[{\rm (D3)}]  $\Ad(\Delta)$ is definable over $\mathbb{A}$, and
\item[{\rm (D4)}]  $\Ad(\Gamma)$ is definable over $\mathbb{A}$.
\end{itemize}
By Proposition \ref{prop:thm2}, (D2) implies (D3).   
By Proposition \ref{prop:thm3}, (D3) implies (D4). 
Then by Proposition \ref{prop:thm2}, Lemma \ref{lem:D2} and (D4) imply (D1).
\end{proof}

Also, we remark that (D1)--(D4) are equivalent under the assumptions.


\subsection{The definability for groups generated by reflections}


\begin{theorem}\label{thm:main2}
Suppose $n\geq 2$ is an  integer. and let $\mathbb{F}$ be an algebraically closed field of characteristic $0$. 
Let $\Gamma \subset \SL^{\pm}(n+1,\mathbb{F})$ be a $\mathbb{F}$-Zariski-dense subgroup generated by finitely many reflections. 
Let $\Delta$ be any subgroup of finite index in $\Gamma$. Then the following 
are equivalent.
\begin{enumerate}
\item[{\rm (a)}] $\Gamma$ is definable over $\mathbb{A}$.
\item[{\rm (a${}'$)}] $\mathbb{A}$ contains all the simple cyclic products of the Cartan matrix of $\Gamma$.
\item[{\rm (b)}] $\Delta$ is definable over $\mathbb{A}$.

\end{enumerate}
\end{theorem}
\begin{proof}
Since $\Gamma$ is Zariski dense in $\SL^{\pm}(n+1,\mathbb{F})$, it is absolutely irreducible.  
By Propositions  \ref{prop:cyclic} and \ref{prop:tr}, (a) and (a${}'$) are equivalent. 
Since the trace of a reflection element in $\SL^{\pm}(n+1,\mathbb{F})$ is $n-1 \ne 0$, 
the equivalence of (a) and (b) follows from Theorem \ref{thm:maingeneral}. 
\end{proof}

\begin{proof}[The proof of Theorem \ref{thm:mainprop}]
Let $\Gamma$ be as in the premise of Theorem \ref{thm:mainprop}. 
By Theorem \ref{thm:zariski},
the Zariski closure $\overline{\Gamma}$ over $\mathbb{R}$ is 
either $\mathrm{O}^{+}(1,n)$ or $\SL^{\pm}(n+1,\mathbb{R})$. 
The result follows from Theorem \ref{thm:main2} or Proposition \ref{thm5}, i.e., Theorem 5 of Vinberg \cite{Vinberg3}.
\end{proof}

\subsection{The definability for groups generated by finite-order elements}




Let $A_{\theta}$ be the matrix
\[\left(
  \begin{array}{cc}
    \cos \theta & -\sin \theta \\
    \sin \theta & \cos \theta 
   \end{array}
\right).\]  
Let $A_{\theta}^{[d]}$ denote the direct sum of $d$, $d \geq 1$, copies of $A_{\theta}$. 
In this section, we state and prove corollaries of Theorem \ref{thm:main2} in the case when the generators are 
the special type of the following well-known lemma. For example see ~\cite{Koo}.
\begin{lemma}\label{lem:theta}
An element $g$ in $\SL^{\pm}(n+1,\mathbb{R})$ has a finite order if and only if $g$ is conjugate to 
\begin{equation}\label{conjugate}
A_{\theta_{1}}^{[d_{1}]} \oplus \cdot \cdot \cdot \oplus A_{\theta_{r}}^{[d_{r}]} \oplus I_{k_{1}} \oplus -I_{k_{2}},
\end{equation}
where 
\begin{align*}
& k_{1}, k_{2} \geq 0, r \geq 0, d_{1}, \cdot \cdot \cdot , d_{r} \geq 1, 0< \theta_{1} < \cdot \cdot \cdot < \theta_{r} < \pi \\
& k_{1}+k_{2}+2(d_{1}+ \cdots d_{r})=n+1 
\end{align*}
hold, and
each $\theta_{i}= 2 \pi a_{i}/b_{i}$ with $(a_{i}, b_{i})=1$ is a rational multiple of $2 \pi$
Then the order of $g$ is $\textrm{lcm}\{2,b_{1} \cdot \cdot \cdot, b_{r}\}$ 
or $\textrm{lcm}\{b_{1} \cdot \cdot \cdot, b_{r}\}$ according as $k_{2}>0$ or $k_{2}=0$ respectively.
\end{lemma}

The conclusion of the following corollary is also true if $n=2$ and the every order $b_i \ne 3$ by the trace condition of 
Theorem \ref{thm:maingeneral}. 
However, this is unclear for the triangle group of orders $(3, 3, 4)$ considered in \cite{Long}. 
Corollary \ref{cor:dim4} is applicable to the projective group isomorphic to 
the fundamental group of the double orbifold of the Coxeter $3$-orbifold based 
on a complete hyperbolic polytope of all dihedral angles $\pi/3$ or the double orbifold of 
a compact Coxeter $n$-orbifold for $n \geq 4$.  (See Chapter 4 of \cite{Choi3} for the definition of doubling and the computing the fundamental groups.)

\begin{corollary}\label{cor:dim4}
 Let $n \geq 3$. Let $\Gamma \subset \SL(n+1,\mathbb{R})$ be a subgroup generated by finitely many elements. 
 Suppose that we have a presentation of $\Gamma$ so that every generator $g_{i}$  
 of $\Gamma$ associated to the presentation is of finite order $b_{i}$ and $g_{i}$ 
 is conjugate to the type $A_{\theta_{i}} \oplus I_{k_{1}}$.
 Suppose that $b_i> 2$ when $n=3$. 
Let $\mathbb{A}$ be an integrally closed Noetherian ring in $\mathbb{R}$, and let $\mathbb{K}$ be the field of fractions of $\mathbb{A}$.  
Suppose that the Zariski closure of $\Gamma$ over $\mathbb{R}$ is $\SL(n+1,\mathbb{R})$.
Let $\Delta$ be any subgroup of finite index in $\Gamma$. 
Then the following are equivalent.
\begin{enumerate}
\item[{\rm (a)}] $\Gamma$ is definable over $\mathbb{A}$.
\item[{\rm (b)}] $\Delta$ is definable over $\mathbb{A}$.
\end{enumerate}
\end{corollary}
\begin{proof}
We claim that $\mathbb{Q} (\mathrm{tr} \Ad g_{i})=\mathbb{Q} (\mathrm{tr} g_{i})$ for each $i$. 
By assumption we have that $g_{i}$ is conjugate to $A_{\theta_i} \oplus I_{k_{1}}$ where $\theta_i=\frac {2a_{i} \pi} {b_{i}}$ 
with $a_{i}$ a positive integer relatively prime to $b_{i}$. Since some power of $A_{\theta_i} \oplus I_{k_{1}}$ is 
$A_{\frac {2\pi} {b_{i}}} \oplus I_{k_{1}}$ and vice versa, we can assume that $\theta_i=\frac {2\pi} {b_{i}}$.
 We have that \[\mathrm{tr} \Ad g_{i}=(\mathrm{tr} g_{i})^{2}-1 = 4\cos^{2} \frac {2\pi} {b_{i}} + 4k_{1} \cos \frac {2\pi} {b_{i}}+k_{1}^{2}-1.\]
 Hence \[\mathbb{Q} (\mathrm{tr} g_{i}) = \mathbb{Q}(\cos \theta_1)\]
 is at most a degree $2$ extension over \[\mathbb{Q}(\mathrm{tr} \Ad g_{i}) = \mathbb{Q}\left(\cos^2 \frac {2\pi} {b_{i}} + k_{1} \cos \frac {2\pi} {b_{i}}\right).\] 
 
 Suppose that $\mathbb{Q} (\mathrm{tr} \Ad g_{i}) \neq \mathbb{Q} (\mathrm{tr} g_{i})$.  
 Then $\mathbb{Q} (\mathrm{tr} g_{i})$ 
is a degree $2$ extension of  $\mathbb{Q} (\mathrm{tr} \Ad g_{i})$. 
Thus there exists a non trivial automorphism $\rho$ of  $\mathbb{Q} (\mathrm{tr} g_{i})$  
which is identity restricted on $\mathbb{Q} (\mathrm{tr} \Ad g_{i})$. Namely, there exists a nontrivial Galois automorphism
\[\rho \in \mathrm{Aut}(\mathbb{Q}(\cos \theta_i)/\mathbb{Q}(\cos^{2}\theta_i+k_{1} \cos \theta_i)).\]  
Since $\cos \theta_i$ and $-k_{1}-\cos \theta_i$ are two roots of the polynomial
\begin{displaymath}
x^{2}+k_1 x-(\cos^{2}\theta_i+ k_1 \cos \theta_i),
\end{displaymath}
we have that $\rho(\cos \theta_i)=-k_{1}-\cos \theta_i$. However since $\mathbb{Q}(\cos \theta_i)=\mathbb{Q}(\zeta_{l}+\zeta_{l}^{-1})$ 
where $\zeta_{l}=e^{i \theta_i}$, we have that $\rho(\cos \theta_i)=\cos m\theta_i$ for some positive integer $m$
by the cyclotomic field theory over $\mathbb{Q}$.

Note that $n \geq 4$ implies $k_{1} \geq 3$ and
\begin{equation}
|\cos m \theta_i| \leq 1 < |-k_{1}-\cos \theta_i|. 
\end{equation}
Now we have a contradiction. 
Also, if $n=3$, then $k_1 =2$. Since $\theta_1 < \pi$, we obtain the same contradiction,
and this proves the claim. 

The claim implies that 
 $\mathrm{tr} g_{i}$ is in $\mathbb{K}$ since $\mathbb{Q} (\mathrm{tr} \Ad \Gamma)$ is the smallest field of definition of $\Ad \Gamma$ 
 and the smallest field of definition of $\Ad \Gamma$ is contained in the smallest field of definition of  $\Gamma$ by ~\cite{Vinberg3}. 
 Note that $n \geq 4$ or $n=3$ with $b_i \ne 2$ for all $i$ implies that $\mathrm{tr} g_{i}$ is never $0$.
Hence the corollary follows from Theorem \ref{thm:maingeneral}. 
\end{proof}

 \section{Examples}\label{sec:examples}

\subsection{A preliminary lemma}
 In this section, for a few $2$- , $3$- and $4$-dimensional Coxeter orbifold $\hat{P}$,
 we find all or some conjugacy classes of irreducible dividing polyhedral projective reflection groups that are definable over $\mathbb{Z}$
 or some other ring of integers. In particular, groups such as $\SL(n+1, \mathbb{Z})$ for $n=2,3,4$ contains the fundamental groups of hyperbolic $n$-manifolds. 
 
 The following lemma of Vinberg ~\cite{Vinberg2} will be useful in this section.
 
 \begin{lemma}\label{lem:distinguish}
 Let $\hat{P}$ be a compact $n$-dimensional Coxeter orbifold. Suppose that $\Gamma_{1}$ and $\Gamma_{2}$ are convex projective reflection groups
 and that $\Omega_{1}$ and $\Omega_{2}$ are irreducible properly convex open domains that the groups divide respectively. Then the following statements are equivalent.
 \begin{enumerate}
\item[{\rm (i)}] $\Omega_{1}/\Gamma_{1}$ and $\Omega_{2}/\Gamma_{2}$ are  projectively diffeomorphic.
\item[{\rm (ii)}] A Cartan matrix $A$ of $\Gamma_{1}$ and a Cartan matrix $B$ of $\Gamma_{2}$ are equivalent.
Namely $A=DBD^{-1}$ for a diagonal matrix $D$ having positive diagonal elements.
\item[{\rm (iii)}] The cyclic products of $A$ and the cyclic products of $B$ are identical.
\end{enumerate}
 \end{lemma}
 \begin{proof}
 The equivalence of (ii) and (iii) follows from Proposition 16 of ~\cite{Vinberg2}. 
The fact that (i) implies (ii) is from the ambiguity of choices in equation \ref{eqn:amb}.
 Since $\Gamma_{1}$ and $\Gamma_{2}$ 
 divide irreducible properly convex open domains, the Vey irreducibility theorem 
 (Theorem 5.1 of ~\cite{Benoistsurvey}) implies that $\Gamma_{1}$ and $\Gamma_{2}$ are irreducible.  
Hence, the Cartan matrix is irreducible. 
 By Corollary 1 of ~\cite{Vinberg2} implies the equivalence of (i) and (ii). 
 \end{proof}


Suppose that $Q$ is a  compact hyperbolic orbifold.
We recall that every finite index subgroup of the fundamental group $\pi_{1}(Q)$ has a trivial center.
 Let $\Gamma$ in $\SL^{\pm}(n+1,\mathbb{R})$ be a group isomorphic to $\pi_{1}(Q)$ dividing a properly convex open domain $\Omega$ in $\mathbb{S}^{n}$.
 Then Corollary 2.13 of Benoist ~\cite{Benoist3} implies that such a group
 $\Gamma$ in $\SL^{\pm}(n+1,\mathbb{R})$ is strongly irreducible so that $\Omega$ has to be irreducible.
 By Theorem 1.1 of ~\cite{Benoist1}, $\Omega$ moreover has to be strictly convex.
 Hence for a compact hyperbolic Coxeter orbifold $\hat{P}$,
 we can apply Theorem \ref{thm:mainprop} to find out whether or not 
 there exists an irreducible dividing polyhedral projective reflection group isomorphic to $\pi_1(\hat{P})$
 that is definable over $\mathbb{Z}$.
 
 Note that $4\cos^{2}(\frac{\pi}{n})$ is an integer if and only if $n \in \{2,3,4,6\}$. 
 By the condition (L2) of Theorem \ref{reflection} and by Theorem \ref{thm:mainprop}, a necessary condition 
 that for an irreducible dividing polyhedral projective reflection group
 to be definable over $\mathbb{Z}$ is that the edge orders are in $\{2,3,4,6\}$.

\subsection{Orbifolds based on simplices}
The \emph{$(p,q,r)$-triangle} is defined to be the triangle with vertices of orders $p, q,$ and $r$. 
The corresponding Coxeter orbifold is said to be a {\em $(p, q, r)$-triangular Coxeter orbifold}.
Let $\hat{P}$ be a $(p,q,r)$-triangular Coxeter orbifold.
It is an elementary fact that $\hat{P}$ admits a hyperbolic structure if and only if  $\frac{1}{p}+\frac{1}{q}+\frac{1}{r}$ is less than $1$.

\begin{proposition}\label{prop:triangle}
Suppose that $\frac{1}{p}+\frac{1}{q}+\frac{1}{r}<1$. Let $\hat{P}$ be a $(p,q,r)$-triangular Coxeter orbifold.
The following list shows the number of conjugacy classes of  irreducible dividing polyhedral projective reflection groups isomorphic to $\pi_1(\hat{P})$ in
$\SL^{\pm}(3,\mathbb{R})$ that are definable over $\mathbb{Z}$.
 \begin{displaymath}
\begin{array}{|c|c||c|c||c|c|}
\hline
\mathrm{triple} & N& \mathrm{triple} & N& \mathrm{triple} & N\\
\hline
(2,4,6) & 1 & (3,4,4) & 3 & (4,4,6) & 6\\
\hline
(3,3,4) & 2 &(3,4,6) & 4 &  (4,6,6) & 5\\
\hline
(3,3,6) & 2 &(4,4,4) & 4 &  (6,6,6) & 4\\
\hline
(2,6,6) & 1 &(3,6,6) & 3 &          &   \\
\hline

\end{array}
\end{displaymath}
Here $N$ is the number of conjugacy classes of  such groups isomorphic to
$\pi_{1}\hat{P}$ definable over $\mathbb{Z}$.
\end{proposition}

\begin{proof}
The list includes every $(p, q, r)$, $p, q, r \in \{2, 3, 4, 6\}$ and $\frac{1}{p}+\frac{1}{q}+\frac{1}{r}<1$.
Let
\begin{displaymath}
\left(
  \begin{array}{ccc}
    2 & c_{12} & c_{13} \\
    c_{21} & 2 & c_{23} \\
    c_{31} & c_{32} & 2 \\
  \end{array}
\right)
\end{displaymath}
be the Cartan matrix of the projective reflection group associated with a triangle orbifold $\hat{P}$ indicated in the above list. 
Suppose that $c_{ij}$ satisfies the conditions (L1) and (L2). 
Then the determinant necessarily equals 
\[ 8 - 8\left( \cos^2 \frac{\pi}{p} +  \cos^2 \frac{\pi}{q} +  \cos^2 \frac{\pi}{r}\right) + c_{12}c_{23}c_{31}  + c_{21} c_{13} c_{32}\]
by our conditions (L1) and (L2).
We obtain that the last two terms are $\leq 0$ by the sign condition and 
the first two terms give a negative number in all our cases above. 
Hence the matrix is of full rank necessarily if (L1) and (L2) are satisfied for some set of $c_{ij}$s. 
For every $(p, q, r)$, the Cartan matrix is irreducible also since there no two $2$s. 

We give a proof when $(p,q,r)=(3,3,4)$. The other cases are completely analogous.
By the condition (L2) of Theorem \ref{reflection}, 
the values of 
$c_{12}c_{21}$, $c_{13}c_{31}$, and $c_{23}c_{32}$ are necessarily $1$, $1$, and $2$ respectively.
Then since the product of $c_{12}c_{21}$, $c_{13}c_{31}$, and $c_{23}c_{32}$ 
are equal to the product of $c_{12}c_{23}c_{31}$ and $c_{13}c_{32}c_{21}$,
we can have only two possible different integer tuples values $(1,1,2,-1,-2), (1,1,2,-2,-1)$ of simple cyclic products
\[(c_{12}c_{21}, c_{13}c_{31}, c_{23}c_{32}, c_{12}c_{23}c_{31}, c_{13}c_{32}c_{21}).\]
Conversely, we certainly find a list of $c_{ij}$ giving each of the two collections of numbers. 
Given such a list of $c_{ij}$, we obtain our properly convex projective orbifold by Corollary 1 of \cite{Vinberg2}
since the determinant is negative. 
Hence we have two conjugacy classes of  irreducible dividing polyhedral projective reflection groups isomorphic to $\hat P$ that are definable over $\mathbb{Z}$
by Lemma \ref{lem:distinguish}. 

We write the matrices in these cases: 
\[\left(
  \begin{array}{ccc}
    2 & -1 & -1 \\
    -1 & 2 & -1 \\
    -1 & -2 & 2 \\
  \end{array}
\right),
\left(
  \begin{array}{ccc}
    2 & -1 & -1 \\
    -1 & 2 & -2 \\
    -1 & -1 & 2 \\
  \end{array}
\right).\]
 
\end{proof}

\begin{remark}
Let $\hat P$ be the orbifold based on triangle with orders $(3, 3, 4)$. 
 The fundamental group has the following Coxeter group presentation,
 \begin{displaymath}
 \langle s_{1}, s_{2}, s_{3}| s_{1}^{2}=s_{2}^{2}=s_{3}^{2}=(s_{1}s_{2})^{3}=(s_{2}s_{3})^{3}=(s_{3}s_{1})^{4}=1 \rangle.
 \end{displaymath}
 Theorem \ref{thm:mainprop} implies that there are only two conjugacy classes 
of the irreducible convex dividing polyhedral projective reflection groups denoted by $\Delta$ which satisfy the following properties.
\begin{itemize}
\item $\Delta$ is isomorphic to the fundamental group.
\item $\Delta$ is definable over $\mathbb{Z}$.
\end{itemize}

We contrast this with the below. 
 Let $\tilde P$ be the double orbifold of $\hat P$ with the fundamental group  
$H$ with the following presentation,
  \begin{displaymath}
 \langle s_{1}, s_{2}| s_{1}^{3}=s_{2}^{3}=(s_{1}s_{2})^{4} =1\rangle.
 \end{displaymath}
Long, Reid, and Thistlethwaite \cite{Long} showed that 
there are infinitely many conjugacy classes of irreducible convex dividing polyhedral projective groups isomorphic to $H$
definable over $\mathbb{Z}$. 

We note that the result was used to proving that $\SL(3, \mathbb{Z})$ contains the fundamental group 
of orientable closed surface of genus $\geq 2$. 
\end{remark}

Sometimes it is possible to find a 1-parameter family of Cartan matrices such that the set of associated irreducible 
dividing polyhedral projective reflection groups contain all which are definable over $\mathbb{Z}$. For example, let
\begin{displaymath}
\left(
  \begin{array}{ccc}
    2 & -3/2 & -1 \\
    -2 & 2 & -2t \\
    -2 & -1/(2t) & 2 \\
  \end{array}
\right)
\end{displaymath}
be the matrix with $t>0$. 
The determinant equals $-6t - 1/t -4 < 0$. 
Then this $1$-parameter family of matrices satisfies the condition to be Cartan matrices for 
the $(3,4,6)$ triangular irreducible  dividing polyhedral projective reflection groups. 
This family contains all of four groups which are definable over $\mathbb{Z}$; 
these are the cases when $t$ is equal to $\frac{1}{6}$, $\frac{1}{3}$, $\frac{1}{2}$, and $1$.  
Figures 1 -- 4 indicate the images of the convex domains which the four groups divide. 
These figures were drawn from TrianglegroupProj2.nb (see \cite{Choi4}).

\begin{figure}[ht]
\begin{minipage}[b]{0.45\linewidth}
\centering
\includegraphics[scale=0.7]{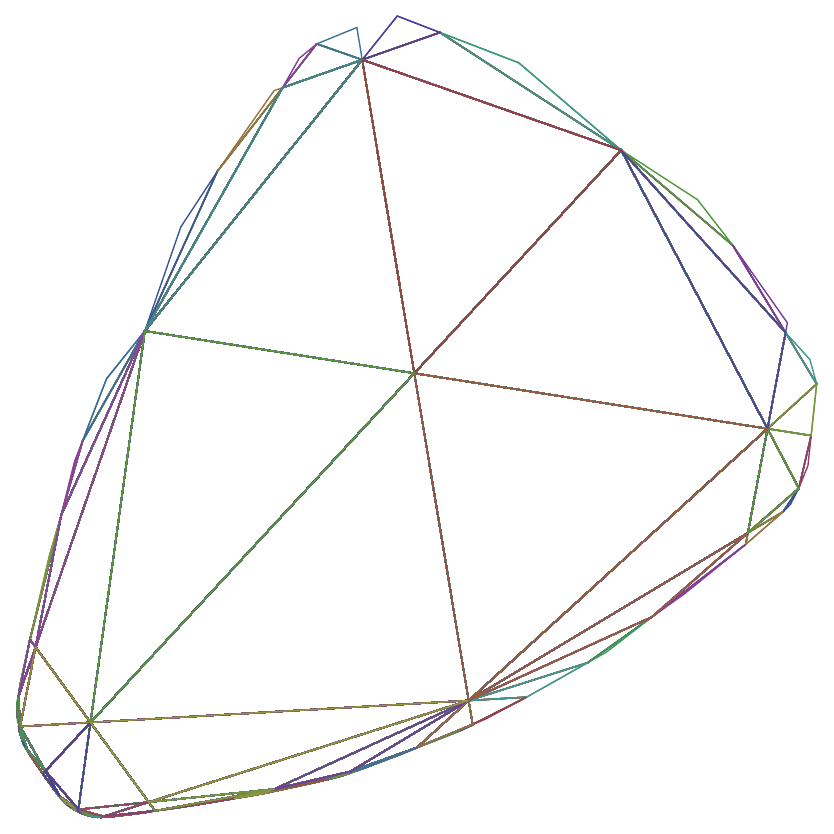}
\caption{$t=\frac{1}{6}$}
\label{fig:figure1}
\end{minipage}
\hspace{0.5cm}
\begin{minipage}[b]{0.45\linewidth}
\centering
\includegraphics[scale=0.7]{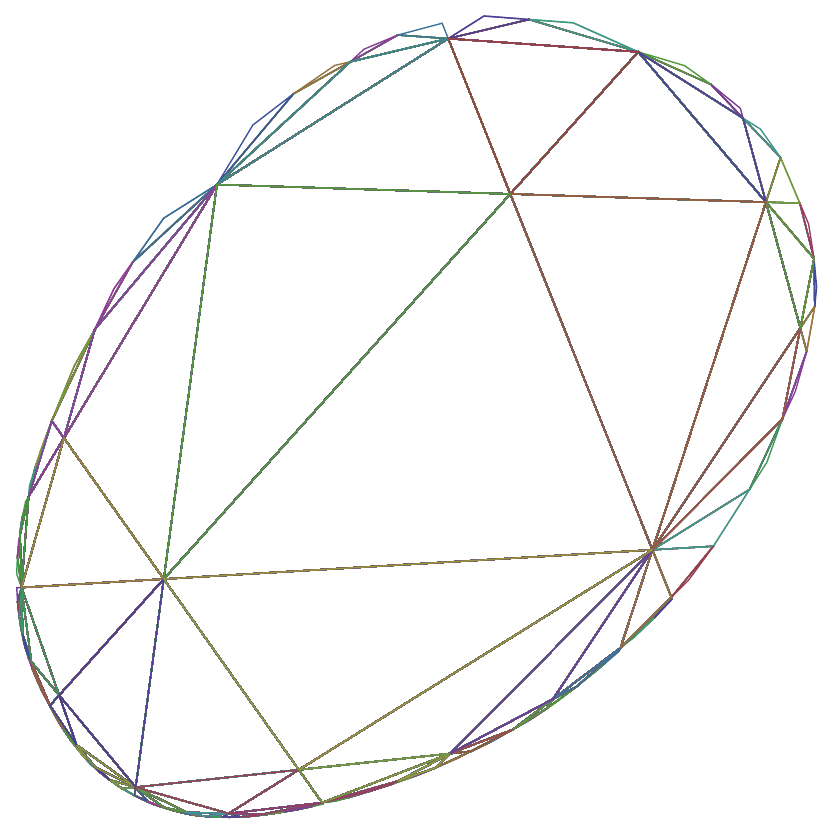}
\caption{$t=\frac{1}{3}$}
\label{fig:figure2}
\end{minipage}
\hspace{0.5cm}
\begin{minipage}[b]{0.45\linewidth}
\centering
\includegraphics[scale=0.7]{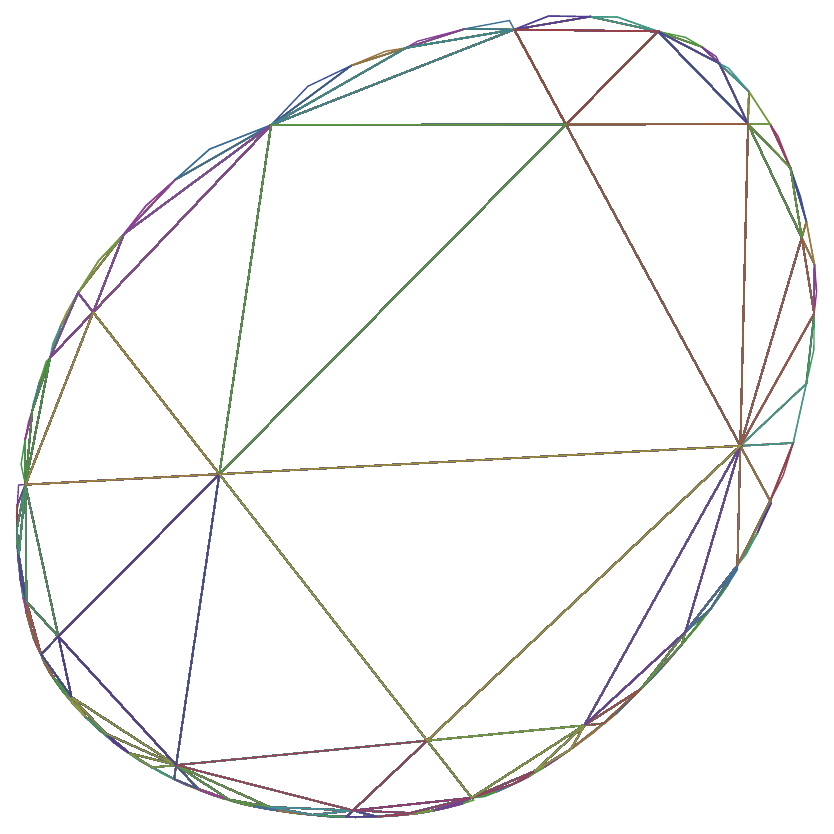}
\caption{ $t=\frac{1}{2}$}
\label{fig:figure3}
\end{minipage}
\hspace{0.5cm}
\begin{minipage}[b]{0.45\linewidth}
\centering
\includegraphics[scale=0.7]{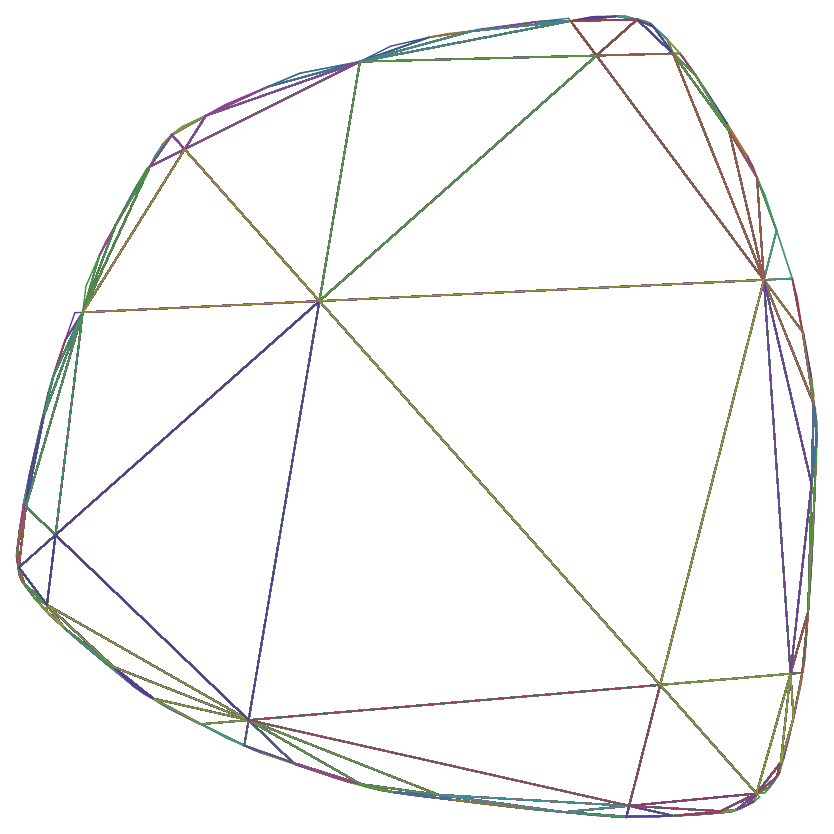}
\caption{$t=1$}
\label{fig:figure4}
\end{minipage}
\end{figure}


\begin{proposition}
Suppose that $\frac{1}{2}+\frac{1}{q}+\frac{1}{r}<1$. Let $\hat{P}$ be a $(2,q,r)$-triangular Coxeter orbifold. 
Then there exists a unique conjugacy class of  
irreducible dividing polyhedral projective reflection groups isomorphic to $\pi_1(\hat{P})$ in $\SL^{\pm}(3,\mathbb{R})$ 
and its smallest ring of definition is $\mathbb{Z}[4\cos^{2}(\frac{\pi}{q}), 4\cos^{2}(\frac{\pi}{r})]$.
\end{proposition}
\begin{proof}
First, there exists such an orbifold by constructing a triangle with angles $\pi/p, \pi/q, \pi/r$ by hyperbolic geometry.
Let
\begin{displaymath}
\left(
  \begin{array}{ccc}
    2 & c_{12} & c_{13} \\
    c_{21} & 2 & c_{23} \\
    c_{31} & c_{32} & 2 \\
  \end{array}
\right)
\end{displaymath}
be the Cartan matrix of such a group which satisfies the equations 
\[c_{12}c_{21}= 4\cos^{2}\frac{\pi}{2} = 0, c_{13}c_{31}= 4\cos^{2}\frac{\pi}{q} \hbox{ and } c_{23}c_{32}=4\cos^{2}\frac{\pi}{r}.\] 
Then the tuple of simple cyclic product 
\[(c_{12}c_{21}, c_{13}c_{31}, c_{23}c_{32}, c_{12}c_{23}c_{31}, c_{13}c_{32}c_{21})\] has to be 
\[(0, 4\cos^{2}\left(\frac{\pi}{q}\right), 4\cos^{2}\left(\frac{\pi}{r}\right), 0, 0).\] 
This determines the  irreducible convex dividing  polyhedral projective  reflection group up to conjugation by Lemma \ref{lem:distinguish}.
The group is absolutely irreducible since the Zariski closure over $\mathbb{R}$ is determined by Theorem \ref{thm:zariski}. 
By Propositions \ref{prop:cyclic} and \ref{prop:tr}, 
$\mathbb{Z}[4\cos^{2}(\frac{\pi}{q}), 4\cos^{2}(\frac{\pi}{r})]$ is its smallest ring of definition.  
\end{proof}

In Section 7.2 of Ratcliffe~\cite{Ratcliffe}, $n$-dimensional simplicial Coxeter orbifolds admitting
hyperbolic structures are classified; indeed there is no $n$-dimensional hyperbolic simplicial orbifold for $n> 4$.  

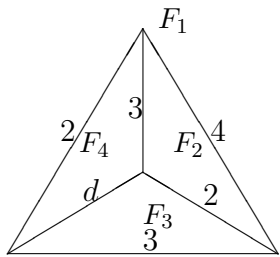
\begin{figure}[h]
\centering
\setlength{\unitlength}{1cm}
\begin{picture}(4,4)
\put(1.8,0.4){\small $F_3$}
\put(0.95,1.35){\small $F_4$}
\put(2.2,1.35){\small $F_2$}
\put(2,3){\small $F_1$}
\put(0,0){\line(1,0){3.6}}
\put(1.8,0.05){$3$}
\put(0,0){\line(3,5){1.8}}
\put(0.7,1.5){$2$}
\put(3.6,0){\line(-3,5){1.8}}
\put(2.7,1.5){$4$}
\put(0,0){\line(5,3){1.8}}
\put(1,0.7){$d$}
\put(3.6,0){\line(-5,3){1.8}}
\put(1.6,1.8){$3$}
\put(1.8,3){\line(0,-1){1.93}}
\put(2.6,0.65){$2$}
\end{picture}
\caption{The compact Coxeter tetrahedron $d=3$ or $4$ }\label{fig:tetrahedron}
\end{figure}

The $n=2$ case is studied in Proposition \ref{prop:triangle}.  
For $n=3$, there are nine simplicial Coxeter orbifolds admitting hyperbolic structures. 
Among them, only two tetrahedral Coxeter orbifolds have edge orders in $\{2, 3, 4, 6\}$. 
The two Coxeter tetrahedrons are described in Figure \ref{fig:tetrahedron}.
For $n=4$, there are only $5$ simplicial Coxeter orbifolds admitting hyperbolic structures . Among them there is only one Coxeter orbifold 
having  edge orders in $\{2, 3, 4, 6\}$. The Coxeter graph of this 4-dimensional simplex is described in Figure \ref{fig:simplex}. 

The following two propositions complete the classification of simplicial orbifolds that admit hyperbolic structures 
whose holonomy groups of convex real projective structures are definable over $\mathbb{Z}$.
(We believe that Goldman has done this in his unpublished notes near 2000.) 

\begin{proposition}\label{prop:tetrahedron}
Let $\hat{P}$ be one of tetrahedral Coxeter orbifolds  described in Figure \ref{fig:tetrahedron}.
If $d=3$, then the number of conjugacy classes of the  irreducible dividing
polyhedral projective  reflection groups isomorphic to $\pi_1(\hat{P})$ in $\SL^{\pm}(4,\mathbb{R})$ that are definable over $\mathbb{Z}$ is $2$.
If $d=4$, then the number of conjugacy classes of the  irreducible dividing polyhedral projective 
reflection groups isomorphic to $\pi_1(\hat{P})$ in $\SL^{\pm}(4,\mathbb{R})$ that are definable over $\mathbb{Z}$ is $3$.
\end{proposition}

\begin{proof}
Let
\begin{displaymath}
\left(
  \begin{array}{cccc}
    2 & c_{12} & c_{13} & c_{14}\\
    c_{21} & 2 & c_{23} & c_{24}\\
    c_{31} & c_{32} & 2 & c_{34}\\
    c_{41} & c_{42} & c_{43}& 2\\
  \end{array}
\right)
\end{displaymath}
be the Cartan matrix of the projective reflection group. 
By the condition (L2) of Theorem \ref{reflection}, 
values of 
\[c_{12}c_{21}, c_{13}c_{31}, c_{14}, c_{41}, c_{23},  c_{32}, c_{24}c_{42}, c_{34}c_{43}\] 
are necessarily \[2, 1, 0, 0, 0, 0, 1, \hbox{ and } 4\cos^{2}\frac{\pi}{d}\] 
respectively. 
Now it is easy to check that every simple cyclic product of length $3$ is $0$, and only nonzero simple cyclic products of length $4$ are 
\[c_{12}c_{24}c_{43}c_{31} \hbox{ and } c_{13}c_{34}c_{42}c_{21}.\]

Finally since the product of $c_{12}c_{21}$, $c_{13}c_{31}$, $c_{24}c_{42}$ and $c_{34}c_{43}$ 
is equal to the product of $c_{12}c_{24}c_{43}c_{31}$ and $c_{13}c_{34}c_{42}c_{21}$, if $d=3$, 
we conclude that only possible integer tuples of nonzero simple cyclic products 
\[(c_{12}c_{21}, c_{13}c_{31}, c_{24}c_{42}, c_{34}c_{43}, c_{12}c_{24}c_{43}c_{31}, c_{13}c_{34}c_{42}c_{21})\] 
is one of $(2, 1, 1, 1, 2, 1)$ and $(2, 1, 1, 1, 1, 2)$.
The determinant is 
\[-c_{12}c_{24}c_{43}c_{31} - 3 c_{34}c_{43} < 0.\]
Clearly, the two possibilities are satisfied for some list $c_{ij}$ of integers
and the Cartan matrices are irreducible. 
We have exactly two different equivalence classes of Cartan matrices whose cyclic products are integers
by Corollary 1 of \cite{Vinberg2} and Lemma \ref{lem:distinguish}.  

If $d=4$, we conclude that only possible integer tuples of nonzero simple cyclic products
\[(c_{12}c_{21}, c_{13}c_{31}, c_{24}c_{42}, c_{34}c_{43}, c_{12}c_{24}c_{43}c_{31}, c_{13}c_{34}c_{42}c_{21})\] 
are
\[(2, 1, 1, 2, 4, 1), (2, 1, 1, 2, 1, 4), \hbox{ and } (2, 1, 1, 2, 2, 2).\]
The corresponding determinants are determined to be negative in these cases. 
Each of these are realized by some list $\{c_{ij}\}$ of integers with irreducible Cartan matrices. 
We have three different equivalence classes of Cartan matrices whose cyclic products are integers
by Corollary 1 of \cite{Vinberg2} and Lemma \ref{lem:distinguish}.  

\end{proof}

\begin{proposition}\label{prop:4-simplex}
Let $\hat{P}$ be the 4-dimensional simplicial Coxeter orbifold  whose Coxeter graph is described in Figure \ref{fig:simplex}.
Then the number of conjugacy classes of the  irreducible dividing polyhedral projective 
reflection groups isomorphic to $\pi_1(\hat{P})$ in $\SL^{\pm}(5,\mathbb{R})$ that are definable over $\mathbb{Z}$ is $2$.

\end{proposition}
\begin{figure}[h]
\[\begin{xy}/r11mm/ :
(2,4.8), {\xypolygon5{@{*}}}
,(2,4.2), {*{4}}
\end{xy}\]
\caption{Coxeter graph of the 4-dimensional simplex definable over $\mathbb{Z}$}\label{fig:simplex}
\end{figure}
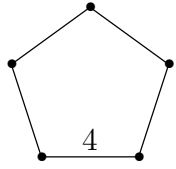
\begin{proof}
Let
\begin{displaymath}
\left(
  \begin{array}{ccccc}
    2 & c_{12} & c_{13} & c_{14} & c_{15}\\
    c_{21} & 2 & c_{23} & c_{24} & c_{25}\\
    c_{31} & c_{32} & 2 & c_{34} & c_{35}\\
    c_{41} & c_{42} & c_{43} & 2 & c_{45}\\
    c_{51} & c_{52} & c_{53} & c_{54} & 2\\
  \end{array}
\right)
\end{displaymath}
be the Cartan matrix of the projective reflection group. By the condition (L2) of Theorem \ref{reflection}, 
$c_{12}c_{21}$,
$c_{13}c_{31}$, $c_{14}c_{41}$, $c_{15}c_{51}$, $c_{23}c_{32}$, $c_{24}c_{42}$, $c_{25}c_{52}$, $c_{34}c_{43}$, $c_{35}c_{53}$, and $c_{45}c_{54}$ are necessarily 
$1$, $0$, $0$, $1$, $1$, $0$, $0$, $2$, $0$, and $1$ respectively. 
Now it is easy to check that every simple cyclic product of length $3$ or $4$ is $0$, and 
only nonzero simple cyclic products of length $5$ are $c_{12}c_{23}c_{34}c_{45}c_{51}$ and $c_{15}c_{54}c_{43}c_{32}c_{21}$.

Finally since the product of $c_{12}c_{21}$, $c_{13}c_{31}$, $c_{24}c_{42}$, $c_{34}c_{43}$ and $c_{45}c_{54}$ 
is equal to the product of $c_{12}c_{23}c_{34}c_{45}c_{51}$ and $c_{15}c_{54}c_{43}c_{32}c_{21}$, 
we conclude that only possible integer tuples of nonzero simple cyclic products 
\[(c_{12}c_{21}, c_{23}c_{32}, c_{34}c_{43}, c_{45}c_{54}, c_{12}c_{23}c_{34}c_{45}c_{51}, c_{15}c_{54}c_{43}c_{32}c_{21})\] equals
\[(1, 1, 1, 2, 1, -1, -2) \hbox{ or } (1, 1, 1, 2, 1, -2, -1).\] 
The two possibilities are realizable by some values $\{c_{ij}\}$ forming irreducible Cartan matrices. 
The corresponding determinants are always $1$ necessarily. 
We have two different equivalence classes of Cartan matrices whose cyclic products are integers 
by Corollary 1 of \cite{Vinberg2} and Lemma \ref{lem:distinguish}.  
\end{proof}


\begin{theorem}\label{thm:infinite}
 Let $\hat{P}$ be a compact $n$-dimensional simplicial Coxeter  orbifold, $n \geq 2$, admitting a hyperbolic structure whose Coxeter graph is not a tree.
 Let $\Gamma$ be a reflection holonomy group of a properly convex projective structure on $\hat{P}$.
 Suppose that $k$ is a real number field that is also a field of definition of $\Gamma$.
 Suppose that $k \neq \mathbb{Q}$.
 Then there are infinitely many conjugacy classes of  irreducible dividing polyhedral projective  reflection groups
 isomorphic to $\pi_1(\hat{P})$ that are definable over the ring of integers $O_{k}$.
\end{theorem}
\begin{proof}
By the classification of hyperbolic simplicial orbifolds (Figure 7.2.9. of ~\cite{Ratcliffe}),
if the Coxeter graph of a hyperbolic simplex is not a tree, then  it is a polygon with $n+1$ vertices
denoted by $v_{i_1}, \dots, v_{i_{n+1}}$. 
Since the Coxeter graph is a polygon, there are exactly two nonzero
cyclic products $c_{i_{1}i_{2}}c_{i_{2}i_{3}}...c_{i_{n}i_{n+1}}c_{i_{n+1}i_{1}}$ 
and $c_{i_{1}i_{n+1}}c_{i_{n+1}i_{n}}...c_{i_{2}i_{1}}$ of length $n+1$. 
Also the Cartan matrix is indecomposable since none of the factors $c_{i_{k}i_{k+1}}$, $c_{i_{k+1}i_{k}}$, 
$c_{i_{n+1}i_{1}}$ and $c_{i_{1}i_{n+1}}$ is $0$.

 We choose 
 $c_{i_{1}i_{2}}$ to be a negative unit $-u$ in $O_{k}$ and $c_{i_{2}i_{1}}$ to
 be $-4u^{-1}\cos^{2}\frac{\pi}{d_{i_{1}i_{2}}}$ where $d_{i_{1}i_{2}}$ is the order corresponding to sides labeled $i_1$ and $i_2$. 
 For other entries, if $j<k$ 
 we choose $c_{i_{j}i_{k}}$ to be $-1$, and  $c_{i_{k}i_{j}}$ to be $-4\cos^{2}\frac{\pi}{d_{i_{j}i_{k}}}$
 provided $|k-j| < 1$, or else they are to be $0$. 
 Then we observe that all the simple cyclic products are
 elements of $O_{k}$. In particular, we obtain 
 \[c_{i_{1}i_{2}}c_{i_{2}i_{3}}...c_{i_{n}i_{n+1}}c_{i_{n+1}i_{1}} =  
 (-1)^{n+1}4u\cos^2\frac{\pi}{d_{i_{1}i_{n+1}}}.\]
 Since there are infinitely many units in $O_{k}$, the cyclic product can take infinitely many values 
 by the choice of $u$. 
 
We can show easily that the determinant of the Cartan matrix  
has only nonzero terms equal to $ \pm 2^{n-2} c_{i_k, i_{k+1}}c_{i_{k+1}, i_k}$ for $k=1, \dots, n$ 
 and the two maximal cyclic terms. 
 \begin{align*}
&\pm c_{i_{1}i_{2}}c_{i_{2}i_{3}}...c_{i_{n}i_{n+1}}c_{i_{n+1}i_{1}} = \pm 4u\cos^{2}\frac{\pi}{d_{i_{j}i_{k}}}\hbox{ and } \\
&\pm c_{i_{n+1}i_{n}}c_{i_{n}i_{n-1}}...c_{i_{2}i_{1}}c_{i_{1}i_{n+1}} = \pm u^{-1}C, C \ne 0. 
\end{align*} 
 For infinitely many $u$, the determinant value of the Cartan matrix is not $0$.
By Corollary 1 of \cite{Vinberg2} and Lemma \ref{lem:distinguish}, the result follows. 
 
\end{proof}

\subsection{Cubical orbifolds}

\begin{proposition}
  Let $\hat{P}$ be the Coxeter orbifold denoted ``cu21'' in ~\cite{CHL} that is also described in Figure \ref{fig:cu21}.
  Then there are at least $3$ conjugacy classes of  irreducible dividing  polyhedral projective  reflection
  holonomy groups in $\SL^{\pm}(4,\mathbb{R})$ that are definable over $\mathbb{Z}$.

 \end{proposition}
  \begin{figure}[h]
\centering
\setlength{\unitlength}{1cm}
\begin{picture}(3,3.5)
\put(1.3,1.35){\large $F_2$}
\put(0.2,1.35){\large $F_5$}
\put(1.3,0.25){\large $F_6$}
\put(2.3,1.35){\large $F_3$}
\put(1.3,2.45){\large $F_1$}
\put(3.4,1.35){\large $F_4$}
    \put(1.45,2.05){$2$}
    \put(1,2){\line(1,0){1}}
    \put(2.05,1.4){$3$}
    \put(2,2){\line(0,-1){1}}
    \put(1.45,0.7){$2$}
    \put(2,1){\line(-1,0){1}}
    \put(0.8,1.4){$2$}
    \put(1,1){\line(0,1){1}}
    \put(0.6,2.45){$3$}
    \put(1,2){\line(-1,1){1}}
    \put(2.25,2.45){$2$}
    \put(2,2){\line(1,1){1}}
    \put(2.25,0.35){$2$}
    \put(2,1){\line(1,-1){1}}
    \put(0.6,0.35){$3$}
    \put(1,1){\line(-1,-1){1}}
    \put(-0.2,1.4){$2$}
    \put(0,0){\line(0,1){3}}
    \put(1.45,3.05){$3$}
    \put(0,3){\line(1,0){3}}
    \put(3.05,1.4){$2$}
    \put(3,3){\line(0,-1){3}}
    \put(1.45,-0.3){$3$}
    \put(3,0){\line(-1,0){3}}
\end{picture}
\caption{\textsf{cu21}}\label{fig:cu21}
\end{figure}
 \begin{proof}
 Gye-Seon Lee has computed that the dimension of infinitesimal restricted deformation space of convex
 real projective structures on $\hat{P}$ is $1$ (see \cite{Lee}).
 This gives a one-parameter family of  irreducible dividing  polyhedral projective  reflection holonomy groups.
We find in ~\cite{KChoi} that the corresponding Cartan matrices are
\begin{displaymath}
  \left(
  \begin{array}{cccccc}
    2 & 0 & 0 & -\frac{2}{2+\sqrt{5}t}& -\frac{2}{2+\sqrt{5}t} & -3+\frac{\sqrt{5}t}{2+2\sqrt{5}t}\\
    0 & 2 & -1 & -\sqrt{6} & 0 & 0\\
    0 & -1& 2 & 0 & -\sqrt{6}& 0\\
    -1-\frac{1}{2}\sqrt{5}t & -\sqrt{6} & 0 & 2 & 0 & -\frac{2+\sqrt{5}t}{2+2\sqrt{5}t} \\
    -1-\frac{1}{2}\sqrt{5}t & 0 & -\sqrt{6} & 0 & 2 & -\frac{2+\sqrt{5}t}{2+2\sqrt{5}t} \\
    -3-\frac{1}{2}\sqrt{5}t & 0 & 0 & -\frac{2(1+\sqrt{5}t)}{2+\sqrt{5}t} & -\frac{2(1+\sqrt{5}t)}{2+\sqrt{5}t} & 2 \\
  \end{array}
  \right)
\end{displaymath}
 where $t$ is a real parameter. This is an irreducible Cartan matrix. 
 We have found that in order for every simple cyclic product to be an integer, 
 only possible values of $t$ are $0$, $-\frac{4}{5\sqrt{5}}$ and  $\frac{4}{\sqrt{5}}$, and these give 
 $3$ different integer tuples of simple cyclic products. Hence there are at least $3$
 conjugacy classes of  irreducible dividing polyhedral projective  reflection groups isomorphic to $\pi_1(\hat{P})$ which are definable over $\mathbb{Z}$
 by Corollary 1 of \cite{Vinberg2} and Lemma \ref{lem:distinguish}. 
 \end{proof}

 \begin{remark}
 The following equalities give the integral cartan matrices that are equivalent to 
 the parameterized Cartan matrix of the cubical projective reflection group when $t$ is $0$, $-\frac{4}{5\sqrt{5}}$ 
 and  $\frac{4}{\sqrt{5}}$ by Lemma \ref{lem:distinguish}. 
 {\small
\begin{displaymath}
  \left(
  \begin{array}{cccccc}
    1 & 0 & 0 & 0& 0 & 0\\
    0 & \frac{2}{\sqrt{6}} & 0 & 0 & 0 & 0\\
    0 & 0& \frac{2}{\sqrt{6}} & 0 & 0& 0\\
    0 & 0 & 0 & 1 & 0 & 0 \\
    0& 0 & 0 & 0 & 1 & 0 \\
    0 & 0 & 0 & 0 & 0 & 1 \\
  \end{array}
  \right)
  \left(
  \begin{array}{cccccc}
    2 & 0 & 0 & -1& -1 & -3\\
    0 & 2 & -1 & -\sqrt{6} & 0 & 0\\
    0 & -1& 2 & 0 & -\sqrt{6}& 0\\
    -1 & -\sqrt{6} & 0 & 2 & 0 & -1 \\
    -1& 0 & -\sqrt{6} & 0 & 2 & -1 \\
    -3 & 0 & 0 & -1 & -1 & 2 \\
  \end{array}
  \right) 
\end{displaymath}
\begin{displaymath}
  \times \left(
  \begin{array}{cccccc}
    1 & 0 & 0 & 0& 0 & 0\\
    0 & \frac{\sqrt{6}}{2} & 0 & 0 & 0 & 0\\
    0 & 0& \frac{\sqrt{6}}{2} & 0 & 0& 0\\
    0 & 0 & 0 & 1 & 0 & 0 \\
    0& 0 & 0 & 0 & 1 & 0 \\
    0 & 0 & 0 & 0 & 0 & 1 \\
  \end{array}
  \right)=
  \left(
  \begin{array}{cccccc}
    2 & 0 & 0 & -1& -1 & -3\\
    0 & 2 & -1 & -2 & 0 & 0\\
    0 & -1& 2 & 0 & -2& 0\\
    -1 & -3 & 0 & 2 & 0 & -1 \\
    -1& 0 & -3 & 0 & 2 & -1 \\
    -3 & 0 & 0 & -1 & -1 & 2 \\
  \end{array}
  \right)
\end{displaymath}

\begin{displaymath}
  \left(
  \begin{array}{cccccc}
    1 & 0 & 0 & 0& 0 & 0\\
    0 & \frac{10}{3\sqrt{6}} & 0 & 0 & 0 & 0\\
    0 & 0& \frac{10}{3\sqrt{6}} & 0 & 0& 0\\
    0 & 0 & 0 & \frac{5}{3} & 0 & 0 \\
    0& 0 & 0 & 0 & \frac{5}{3} & 0 \\
    0 & 0 & 0 & 0 & 0 & 5 \\
  \end{array}
  \right)
  \left(
  \begin{array}{cccccc}
    2 & 0 & 0 & -\frac{5}{3}& -\frac{5}{3} & -5\\
    0 & 2 & -1 & -\sqrt{6} & 0 & 0\\
    0 & -1& 2 & 0 & -\sqrt{6}& 0\\
    -\frac{3}{5} & -\sqrt{6} & 0 & 2 & 0 & -3 \\
    -\frac{3}{5}& 0 & -\sqrt{6} & 0 & 2 & -3 \\
    -\frac{13}{5} & 0 & 0 & -\frac{1}{3} & -\frac{1}{3} & 2 \\
  \end{array}
  \right)
\end{displaymath}
\begin{displaymath}
  \times \left(
  \begin{array}{cccccc}
    1 & 0 & 0 & 0& 0 & 0\\
    0 & \frac{3\sqrt{6}}{10} & 0 & 0 & 0 & 0\\
    0 & 0& \frac{3\sqrt{6}}{10} & 0 & 0& 0\\
    0 & 0 & 0 & \frac{3}{5} & 0 & 0 \\
    0& 0 & 0 & 0 & \frac{3}{5} & 0 \\
    0 & 0 & 0 & 0 & 0 & \frac{1}{5} \\
  \end{array}
  \right)=
  \left(
  \begin{array}{cccccc}
    2 & 0 & 0 & -1& -1 & -1\\
    0 & 2 & -1 & -2 & 0 & 0\\
    0 & -1& 2 & 0 & -2& 0\\
    -1 & -3 & 0 & 2 & 0 & -1 \\
    -1& 0 & -3 & 0 & 2 & -1 \\
    -13 & 0 & 0 & -1 & -1 & 2 \\
  \end{array}
  \right)
\end{displaymath}

\begin{displaymath}
 \left(
  \begin{array}{cccccc}
    1 & 0 & 0 & 0& 0 & 0\\
    0 & \frac{2}{3\sqrt{6}} & 0 & 0 & 0 & 0\\
    0 & 0& \frac{2}{3\sqrt{6}} & 0 & 0& 0\\
    0 & 0 & 0 & \frac{1}{3} & 0 & 0 \\
    0& 0 & 0 & 0 & \frac{1}{3} & 0 \\
    0 & 0 & 0 & 0 & 0 & \frac{1}{5} \\
  \end{array}
  \right)
  \left(
  \begin{array}{cccccc}
    2 & 0 & 0 & -\frac{1}{3}& -\frac{1}{3} & -\frac{13}{5}\\
    0 & 2 & -1 & -\sqrt{6} & 0 & 0\\
    0 & -1& 2 & 0 & -\sqrt{6}& 0\\
    -3 & -\sqrt{6} & 0 & 2 & 0 & -\frac{3}{5} \\
    -3& 0 & -\sqrt{6} & 0 & 2 & -\frac{3}{5} \\
    -5 & 0 & 0 & -\frac{5}{3} & -\frac{5}{3} & 2 \\
  \end{array}
  \right)
\end{displaymath}
\begin{displaymath}
 \times  \left(
  \begin{array}{cccccc}
    1 & 0 & 0 & 0& 0 & 0\\
    0 & \frac{3\sqrt{6}}{2} & 0 & 0 & 0 & 0\\
    0 & 0& \frac{3\sqrt{6}}{2} & 0 & 0& 0\\
    0 & 0 & 0 & 3 & 0 & 0 \\
    0& 0 & 0 & 0 & 3 & 0 \\
    0 & 0 & 0 & 0 & 0 & 5 \\
  \end{array}
  \right)=
  \left(
  \begin{array}{cccccc}
    2 & 0 & 0 & -1& -1 & -13\\
    0 & 2 & -1 & -2 & 0 & 0\\
    0 & -1& 2 & 0 & -2& 0\\
    -1 & -3 & 0 & 2 & 0 & -1 \\
    -1& 0 & -3 & 0 & 2 & -1 \\
    -1 & 0 & 0 & -1 & -1 & 2 \\
  \end{array}
  \right)
\end{displaymath}
}

 \end{remark}

\subsection{The prismatic orbifolds of Benoist}

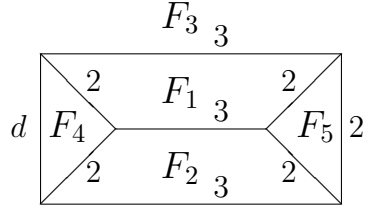
\begin{figure}[h]
\centering
\setlength{\unitlength}{1cm}
\begin{picture}(4,2.7)
\put(0.1,0.9){\large $F_4$}
\put(3.4,0.9){\large $F_5$}
\put(1.6,2.35){\large $F_3$}
\put(1.6,1.35){\large $F_1$}
\put(1.6,0.35){\large $F_2$}
\put(-0.4,0.9){$d$}
\put(0,0){\line(0,1){2}}
\put(0.6,0.3){$2$}
\put(0,0){\line(1,1){1}}
\put(0.6,1.5){$2$}
\put(0,2){\line(1,-1){1}}
\put(2.3,2.1){$3$}
\put(0,2){\line(1,0){4}}
\put(2.3,1.1){$3$}
\put(0,0){\line(1,0){4}}
\put(2.3,0.1){$3$}
\put(4,0){\line(0,1){2}}
\put(3.2,1.5){$2$}
\put(4,0){\line(-1,1){1}}
\put(3.2,0.3){$2$}
\put(4,2){\line(-1,-1){1}}
\put(4.1,0.9){$2$}
\put(1,1){\line(1,0){2}}
\end{picture}
\caption{The prism}\label{fig:triprism}
\end{figure}

So far, all the examples of  irreducible dividing polyhedral projective reflection holonomy  groups we have considered were ones that divide strictly convex domains. Our final
example illustrates that we can apply our main theorem (Theorem \ref{thm:main2}) even when the domain which an  irreducible dividing polyhedral projective reflection  holonomy
group divides is not strictly convex.
\begin{proposition}
 Suppose that $d$ is $3$ or $4$. Let $\hat{P}$ be a triangular prism with dihedral angles described in Figure \ref{fig:triprism}.
 Then the number of  irreducible convex dividing polyhedral projective reflection  holonomy groups in $\SL^{\pm}(4,\mathbb{R})$ which are definable over $\mathbb{Z}$ is $0$.
 \end{proposition}

\begin{proof}
This example is from Benoist ~\cite{Benoist4}. 
Benoist proved that there exists a family of  irreducible dividing polyhedral reflection
groups Zariski-dense in $\SL^{\pm}(4,\mathbb{R})$ whose Cartan matrices are of from
\begin{displaymath}
\left(
  \begin{array}{ccccc}
    2 & -1 & -1 & 0 & 0\\
    -1 & 2 & -t & 0 & 0\\
    -1 & -t^{-1}& 2 & \frac{\mu(1-t)}{t} & 0\\
    0 & 0 & 1-t & 2 & -2\\
    0 & 0 & 0 & -\nu & 2\\
  \end{array}
 \right)
\end{displaymath}
 where $t>1$, $\mu=\frac{4t}{(t-1)^{2}}\cos^{2}(\frac{\pi}{d})$, and $\nu=2+3\mu$.
(See Proposition 4.2 of \cite{Benoist4}.)
 
 We observe that $-t$ and $-t^{-1}$ are both the values of two simple cyclic products. 
 Since it is impossible for both $-t$ and
 $-t^{-1}$ to be integers, Theorem \ref{thm:main2} implies 
 that there is no dividing polyhedral reflection group isomorphic to $\pi_1(\hat{P})$ which is definable over $\mathbb{Z}$.  
 
\end{proof}